\documentclass[11pt]{article}

\usepackage[utf8]{inputenc}
\usepackage[T1]{fontenc}
\usepackage[a4paper, margin=1in]{geometry}
\usepackage{lmodern}
\usepackage{amsmath, amssymb, amsthm}
\usepackage{graphicx}
\usepackage{float}
\usepackage{enumitem}
\usepackage{setspace}
\usepackage{microtype}
\usepackage[most]{tcolorbox}
\usepackage{adjustbox}
\usepackage{booktabs}

\usepackage{tabularx}
\usepackage{hyperref}
\usepackage{pifont}
\newcommand{\xmark}{\ding{55}}  

\hypersetup{
    colorlinks=true,
    linkcolor=blue,
    citecolor=blue,
    urlcolor=blue
}

\theoremstyle{definition} 
\newtheorem{theorem}{Theorem}[section]
\newtheorem{lemma}[theorem]{Lemma}
\newtheorem{corollary}[theorem]{Corollary}
\newtheorem{definition}[theorem]{Definition}
\newtheorem{remark}[theorem]{Remark}

\title{A Finite-State Symbolic Automaton Model for the Collatz Map\\ and Its Convergence Properties}

\author{
  Leonard Ben Aurel Brauer \\
  \small Independent Researcher \\
  \small \texttt{leonard.brauer.research@gmail.com} \\
  \small \href{https://orcid.org/0009-0003-9173-4009}{ORCID: 0009-0003-9173-4009} \\
}

\date{}  

\begin{document}
\maketitle

\textbf{Abstract} 

We present a finite-state, deterministic automaton that emulates the Collatz function through digitwise transitions on base-10 representations. Each digit is represented as a symbolic triplet $(r, p, c)$ encoding its value, the parity of the next digit, and an incoming carry propagated from the lower digit. This yields exactly 60 possible local states. The automaton applies local, parity-aware rules that collectively reconstruct the global arithmetic of the Collatz map. We show that all symbolic trajectories converge in finitely many steps to a unique terminal cycle $(4, 0, 0) \rightarrow (2, 0, 0) \rightarrow (1, 0, 0)$, with all higher digit positions degenerating to the absorbing state $(0, 0, 0)$. This collapse reveals a canonical symbolic normal form of Collatz dynamics.

In parallel, a binary view explains the dynamics as alternating bit-length growth and contraction, aligning with known heuristics for Collatz convergence. This structural perspective is further reinforced by a symbolic drift function and a ranking potential that together explain and formalize the convergence process.


\section{Introduction}

Few problems in elementary mathematics are as deceptively simple and notoriously resistant as the so-called \emph{3\(x + 1\)} problem. Defined by an innocuous piecewise function over the natural numbers, it asks whether repeated application of the map
\[
T(n) = 
\begin{cases}
\frac{n}{2},     & \text{if } n \equiv 0 \pmod{2} \\
3n + 1,          & \text{if } n \equiv 1 \pmod{2}
\end{cases}
\]
will eventually reach the value \(1\) for all \(n \in \mathbb{N}^+\). Despite its apparent simplicity and extensive numerical verification, the conjecture—first posed by Lothar Collatz in 1937—remains unproven \cite{lagarias2003}.
The difficulty lies in the absence of a known invariant or monotonicity property: trajectories exhibit nontrivial recursive fluctuations, and the function's alternating arithmetic structure complicates analytic treatment. In this work, we present a finite-state, symbolic reformulation of the Collatz dynamics. Each base-10 digit of \( n \) is represented as a triplet \((r, p, c)\), where \( r \) is the digit's value, \( p \) the parity of the next higher-order digit, and \( c \) a carry propagated from lower positions. This yields a deterministic automaton with exactly $10 \cdot 2 \cdot 3 = 60$ distinct states. Transitions act locally on digit positions, yet collectively reproduce the global behavior of $T(n)$ via structured, digitwise updates. Each number is thus represented as a vertical stack of symbolic digit states, one per decimal place, where each column evolves independently over time. The symbolic dynamics proceed row-by-row through successive configurations, allowing the Collatz behavior to emerge from purely local interactions.
We prove that all symbolic digit states eventually collapse to the null state \((0,0,0)\), except for the least significant position. Once all higher digit positions are absorbed, the remaining state at position \(s_0\) enters the unique three-node cycle \((4, 0, 0) \rightarrow (2, 0, 0) \rightarrow (1, 0, 0)\), which corresponds to the classical terminal loop \(4 \rightarrow 2 \rightarrow 1\). No other symbolic cycles or divergent paths exist within this structure.

Beyond its theoretical implications, the automaton provides a symbolic and transparent model of integer dynamics. Its local, deterministic transitions offer a natural fit for logic circuits, low-power symbolic processors, or explainable symbolic AI frameworks. The construction opens new perspectives on digit-level number theory and computationally verifiable transformations.

\paragraph{Origin of the Approach.}
This work did not begin with the goal of addressing the Collatz conjecture directly, but grew out of a broader interest in how structured behavior can arise from simple symbolic rules. When the author first encountered the Collatz function—previously unfamiliar—it was immediately clear that its complex dynamics likely stem from its minimal definition, rather than from randomness. Motivated by this, a small tool was developed to visualize and explore Collatz sequences. This quickly revealed patterns at the digit level, pointing toward an underlying local regularity. One key insight was the central role of the least significant digit, \( r_0 \): its value, together with parity and carry information, sufficed to determine the local update behavior. This prompted the design of a symbolic representation, which evolved into the deterministic automaton formalized in this study.

\subsection{Related Work}

Numerous past and ongoing projects in number theory, dynamical systems, and symbolic computation have investigated the behavior of arithmetic iterations resembling the Collatz map, particularly with the aim of proving or disproving convergence under various formulations. An updated overview of historical and recent developments is provided by Lagarias~\cite{lagarias2021overview}, complementing his earlier annotated bibliographies~\cite{lagarias2003, lagarias2006}. Residue class analysis and tree-based trajectory studies emphasize the combinatorial complexity of such functions.
Tao~\cite{tao2019collatz} introduced a probabilistic and ergodic perspective on integer dynamics, studying the behavior of iterative maps under statistical models. This probabilistic viewpoint extends earlier work by Lagarias and Weiss ~\cite{lagarias1992stochastic}, who developed stochastic models using biased random walks to examine the statistical behavior of Collatz trajectories, including growth rates and stopping times.
The present work deviates from these global strategies by proposing a digitwise symbolic emulation using a deterministic finite-state automaton. Each digit of the input number is processed locally, enabling a structural decomposition of the transformation into symbolic transitions with bounded carry propagation.
The automaton construction, along with a symbolic descent function that guarantees convergence, was developed by the author. This fine-grained representation departs from traditional whole-number views and establishes termination via symbolic state space dynamics.
Foundational tools from discrete mathematics and automata theory~\cite{rosen2012discrete, sipser2005automata} underlie the definition of state transitions, carry mechanisms, and formal verification techniques used in this framework.

\section{A Finite-State Symbolic Automaton for the Collatz Function}

We present a symbolic representation of the Collatz function acting on the base-10 digit expansion of natural numbers. This formulation encodes each digit with auxiliary state information, enabling the construction of a deterministic finite-state automaton with 60 states. Each symbolic state encodes the local configuration of a single decimal digit during a Collatz transformation step. Formally, a state is a triplet $(r_i, p_i, c_i)$ where $r_i$ is the $i$-th digit of the input number (in base 10), $p_i = r_{i+1} \bmod 2$ is the parity of the \emph{next} digit and $c_i$ is a carry value originating from the previous digit's computation. These states are processed sequentially from least to most significant digit, enabling the digitwise reconstruction of operations such as $3n + 1$ and $n/2$ using local arithmetic and carry propagation only. The total configuration of a number is thus captured by a sequence of symbolic states indexed by digit position. 

\subsection{Positional Representation and Symbolic State}

Let $n \in \mathbb{N}$ be written as a base-10 expansion \cite{rosen2012discrete}:
\[
n = \sum_{i=0}^L r_i \cdot 10^i, \quad r_i \in \{0,1,\dots,9\}
\]
where $r_0$ is the least significant digit and $L \in \mathbb{N}$ is the index of the most significant digit.

Each digit $r_i$ is associated with:
\begin{itemize}
    \item its value $r_i \in \{0, 1, \dots, 9\}$,
    \item the parity of the next digit: $p_i := r_{i+1} \bmod 2$ for $i < L$; the parity value is undefined (and unused) for $i = L$,
    \item a \emph{carry} $c_i \in \{0,1,2\}$ propagated from the previous digit position, which arises as follows:
    \begin{itemize}
        \item In the \textbf{even case} \( T(n) = \lfloor n/2 \rfloor \), carries arise from parity-dependent digit borrowing.
        \item In the \textbf{odd case} \( T(n) = 3n + 1 \), carries result from digitwise multiplication with overflow.
        \item For the least significant digit \( r_0 \), the initial carry is fixed: \( c_0 := 0 \).
    \end{itemize}

\end{itemize}

We define the symbolic digit state as the triplet $s_i := (r_i, p_i, c_i) \in S$, where the finite state space is given by:
\[
S := \{0,\dots,9\} \times \{0,1\} \times \{0,1,2\}, \quad |S| = 60
\]

From this, a symbolic representation was developed that evolves into a finite-state automaton~\cite{sipser2005automata}, in which each state encodes digit value, parity context, and local carry behavior.

\subsection{Symbolic Transition Rules} \label{sec:transition}

We define $\delta_i : S \to \{0, \dots, 9\}$ as the local digitwise output function. 
It maps each symbolic input state $(r_i, p_i, c_i)$ to the corresponding output digit $r_i'$.
The choice of update rule depends on the global parity $\pi(n) := n \bmod 2$ of the input number:
\[
r_i' := 
\begin{cases}
\delta_{\text{even}}(r_i, p_i, c_i), & \text{if } n \text{ is even}, \\
\delta_{\text{odd}}(r_i, c_i),       & \text{if } n \text{ is odd}.
\end{cases}
\]

Note that $\delta_i$ only computes the updated digit value $r_i'$. The corresponding parity and carry values for subsequent digits are determined separately via recursive rules.

\subsection[Even Case: T(n) = ⌊n/2⌋]{Even Case: $T(n) = \lfloor n/2 \rfloor$}

For even $n$, the Collatz map halves the number. We simulate this operation digit-by-digit using a parity-aware rule.

\subsubsection*{Lookup Table Representation}

The digitwise transformation uses a finite table that maps even digits $r_i \in \{0,2,4,6,8\}$ and parities $p_i \in \{0,1\}$ to a new symbolic digit $r_i'$:

\[
r_i' = \texttt{LOOKUP}\left[\, (r_i - c_i) \bmod 10 \,\right][\,p_i\,]
\]

\begin{table}[h!]
\centering
\caption{Lookup table for even-case digit updates: columns correspond to current digit $r_i$, and rows to the next-digit parity $p_i$.}
\label{tab:lookup-horizontal}
\begin{tabular}{c|ccccc}
$p_i$ & $r_i = 0$ & $r_i = 2$ & $r_i = 4$ & $r_i = 6$ & $r_i = 8$ \\
\hline
0 & 0 & 1 & 2 & 3 & 4 \\
1 & 5 & 6 & 7 & 8 & 9 \\
\end{tabular}
\end{table}
The lookup table defines the transformation for even digits. It maps the adjusted digit \((r_i - c_i) \bmod 10\) and the next-digit parity \(p_i\) to the output digit \(r_i'\).
Each entry specifies the result of halving an even digit $r_i$ under the influence of the parity $p_i$ of the next digit. For example, if $n = 30$, then $r_0 = 0$, $p_0 = 1$, and $c_0 = 0$, leading to $r_0' = \texttt{LOOKUP}[0][1] = 5$.

\subsubsection*{Closed-Form Rule with Carry}

The same mapping can be equivalently expressed as a closed-form update, taking into account the carry $c_i$ propagated from the right:

\[
\delta_{\text{even}}(r_i, p_i, c_i) := \left\lfloor \frac{r_i - c_i}{2} \right\rfloor + 5p_i
\]

This rule matches the lookup behavior and ensures correct adjustment for carry-induced subtractions.

\paragraph{Parity and Carry Propagation:}

\[
p_i := r_{i+1} \bmod 2 \qquad
c_{i+1} := \begin{cases}
1 & \text{if } p_i = 1 \\
0 & \text{otherwise}
\end{cases}
\]

This ensures that each local state $(r_i, p_i, c_i)$ includes the information required to account for parity-induced division shifts, analogous to digit borrowing in manual division.

\subsection{Odd Case: $T(n) = 3n + 1$}

If $n$ is odd, the Collatz map applies $3n + 1$. We emulate this with digitwise multiplication and carry propagation.

\subsubsection*{Digitwise Transition Rule}

 For each digit index $i \geq 0$, define:

\[
\delta_{\text{odd}}(r_i, c_i) := 
\begin{cases}
(3r_i + 1 ) \bmod 10 & \text{if } i = 0 \\
(3r_i + c_i) \bmod 10 & \text{if } i > 0
\end{cases}
\]

\[
c_{i+1} := \begin{cases}
\left\lfloor \frac{3r_i + 1}{10} \right\rfloor & \text{if } i = 0 \\
\left\lfloor \frac{3r_i + c_i}{10} \right\rfloor & \text{if } i > 0
\end{cases}
\]

\subsubsection*{Alternative Representation}

Let $n = \sum_{i=0}^L r_i \cdot 10^i$ be the decimal expansion. Then:
\[
T(n) = 3n + 1 = \left( \sum_{i=0}^{L} 3r_i \cdot 10^i \right) + 1
\]

This expression is evaluated digitwise with carry propagation as per the recursive rules above. The symbolic update faithfully emulates this arithmetic.

\subsection{Unified Digitwise Reconstruction of \boldmath$T(n)$}

Let \( n \in \mathbb{N}^+ \) have decimal expansion \( n = \sum_{i=0}^L r_i \cdot 10^i \), and define:

\[
r_i' := 
\begin{cases}
\delta_{\text{even}}(r_i, p_i, c_i) & \text{if } n \equiv 0 \pmod{2} \\
\delta_{\text{odd}}(r_i, c_i) & \text{if } n \equiv 1 \pmod{2}
\end{cases}
\]

Then the Collatz map is reconstructed symbolically via:

\[
T(n) = \sum_{i=0}^{L'} r_i' \cdot 10^i \quad \text{for some } L' \in \{L - 1, L, L + 1\}
\]
This unified notation allows symbolic digitwise emulation of both cases of the Collatz map.

\subsection{Illustrative Example: $n = 27$}

We illustrate the symbolic digitwise emulation for the odd input \( n = 27 \), where \( T(27) = 3 \cdot 27 + 1 = 82 \).  
The decimal expansion of \( n \) is:
\[
n = 7 \cdot 10^0 + 2 \cdot 10^1 \quad \Rightarrow \quad r_0 = 7,\ r_1 = 2
\]

We now apply the symbolic update rules for the odd case:

\begin{align*}
r_0' &= \delta_{\text{odd}}(7, 0) = (3 \cdot 7 + 1) \bmod 10 = 2 \\
c_1 &= \left\lfloor \frac{3 \cdot 7 + 1}{10} \right\rfloor = 2 \\
r_1' &= \delta_{\text{odd}}(2, 2) = (3 \cdot 2 + 2) \bmod 10 = 8 \\
c_2 &= \left\lfloor \frac{3 \cdot 2 + 2}{10} \right\rfloor = 0
\end{align*}

\begin{table}[H]
\centering
\caption{Symbolic digitwise transition for $n = 27$ under \( T(n) = 3n + 1 \)}
\label{tab:example-27}
\begin{tabular}{c|c|c|c|c|c|c}
$i$ & $r_i$ & $p_i$ & $c_i$ & $r_i'$ & $p_i' = r_i' \bmod 2$ & $c_{i+1}$ \\
\hline
0 & 7 & 0 & 0 & 2 & 0 & 2 \\
1 & 2 & – & 2 & 8 & 0 & 0 \\
\end{tabular}
\end{table}

We obtain the updated digit sequence \( (r_0', r_1') = (2, 8) \), which corresponds to:
\[
T(27) = 2 \cdot 10^0 + 8 \cdot 10^1 = 82
\]

A worked example for an even input (e.g., $n = 32$) is included in Appendix~\ref{appendix:even-example}.

\begin{remark}[Arithmetic Transparency]
The symbolic transition rules emulate classical base-10 arithmetic:

\begin{itemize}
    \item The odd case simulates $T(n) = 3n + 1$ via long multiplication with carry tracking.
    \item The even case performs division by 2 through lookup and parity-based adjustments, analogous to manual division.
\end{itemize}

Hence, the automaton replicates $T(n)$ precisely at the digit level.
\end{remark}

\section{Correctness of Digitwise Transition Rules}

We now confirm that the unified digitwise rule introduced above faithfully emulates the arithmetic behavior of the Collatz function \( T(n) \).

The symbolic update functions \( \delta_{\text{even}} \) and \( \delta_{\text{odd}} \) reproduce the digit-level base-10 representation of \( T(n) \in \mathbb{N} \) exactly, including carry propagation and parity logic. These rules were derived directly from the standard algorithms for division and multiplication in positional notation \cite{knuth1997art}.

\begin{remark}
The correctness of this digitwise formulation implies that the symbolic automaton is not an approximation of Collatz dynamics, but an exact reformulation over a finite-state digit system. This makes it suitable as a rigorous foundation for further formal analysis, including symbolic termination arguments.
\end{remark}

\section{Finite-State Symbolic Convergence (Decimal FSM Model)} \label{sec:convergence}

We now formally prove that the symbolic Collatz automaton terminates in a unique 3-cycle for all $n \in \mathbb{N}^+$. This establishes a fully symbolic convergence proof for the Collatz dynamics, given the numerical correctness of the digitwise emulation (cf. Section~\ref{sec:transition}). The convergence behavior of the Collatz map has been the subject of extensive study, and is considered one of the most notorious open problems in elementary number theory \cite{lagarias2021overview}.

\subsection{Symbolic Encoding and Digitwise Dynamics}

\begin{definition}[Symbolic State Representation]
Let \( n \in \mathbb{N}^+ \) have decimal expansion \( n = \sum_{i=0}^L r_i \cdot 10^i \), where \( r_i \in \{0, \dots, 9\} \).  
We define its symbolic encoding as a sequence of digitwise states
\[
s = [(r_0, p_0, c_0), \dots, (r_L, p_L, c_L)] \in S^{L+1},
\]
where
\begin{itemize}
    \item $p_i := r_{i+1} \bmod 2$ is the parity of the next digit (set to $0$ if $i = L$),
    \item $c_0 := 0$ and subsequent carry values $c_{i+1}$ are determined deterministically by the rules in Section~\ref{sec:transition}.
\end{itemize}
\end{definition}

\begin{definition}[Global Parity]
Given a symbolic state sequence \( s \), we define its global parity by \( \pi(s) := r_0 \bmod 2 \).  
This determines whether the odd or even digitwise update rule applies.
\end{definition}

\begin{theorem}[Correctness of Symbolic Emulation]
\label{thm:encoding}
For every \( n \in \mathbb{N}^+ \), the symbolic encoding \( s \in S^* \) is uniquely defined, and the digitwise update rules $\delta_{\text{even}}$ and $\delta_{\text{odd}}$ reconstruct the next iterate $T(n)$ of the Collatz function via
\[
T(n) = \sum_{i=0}^{L'} r_i' \cdot 10^i, \quad \text{for some } L' \in \{L - 1, L, L + 1\},
\]
where each digit $r_i'$ is computed locally from $(r_i, p_i, c_i)$.
\end{theorem}

\begin{proof}
The base-10 expansion of \( n \) is unique, and both the parity values \( p_i \) and carry values \( c_i \) are deterministically assigned. The digitwise update functions $\delta_i$ compute the output digits $r_i'$ in accordance with standard arithmetic. Therefore, $T(n)$ is reconstructed exactly from the symbolic sequence.
\end{proof}

\begin{remark}
The symbolic encoding map $\mathbb{N}^+ \leftrightarrow S^*$ is bijective.  
Each encoding begins with $c_0 = 0$, and the set of valid initial states is limited to the 20 combinations in $\{0,\dots,9\} \times \{0,1\} \times \{0\}$.
\end{remark}

\subsection{Automaton Structure and Terminal Cycle}

\begin{lemma}[Local Determinism]
The digitwise update rules
\[
\delta_{\text{even}} : S \to \{0, \dots, 9\}, \quad \delta_{\text{odd}} : S \to \{0, \dots, 9\}
\]
are total and deterministic. Carry values \( c_i \in \{0,1,2\} \) are uniquely determined by recursive rules, and thus each symbolic transition is well-defined.
\end{lemma}

\begin{lemma}[Bounded Carry Propagation]
\label{lem:carry-bounds}
For all $s_i \in S$, the carry values $c_i$ remain bounded: $c_i \in \{0, 1, 2\}$.

\noindent\textbf{Odd Case:}  
For worst-case digits \( r_i = 9 \) and \( c_i = 2 \), the carry update (cf. Section~\ref{sec:transition}) yields:
\[
3 \cdot 9 + 2 = 29 \quad \Rightarrow \quad c_{i+1} = \left\lfloor \frac{29}{10} \right\rfloor = 2
\]
At the initial digit \( i = 0 \), with the additional constant term:
\[
3 \cdot 9 + 1 + 0 = 28 \quad \Rightarrow \quad c_1 = \left\lfloor \frac{28}{10} \right\rfloor = 2
\]

\noindent\textbf{Even Case:}  
Carries are computed using the parity of the next digit:
\[
c_i := \begin{cases}
1 & \text{if } r_{i+1} \bmod 2 = 1 \\
0 & \text{otherwise}
\end{cases}
\]

\textbf{Summary:}
In all cases, the carry values are bounded:
\[
c_i \in \{0, 1, 2\} \quad \text{for all } i
\]

\end{lemma}

\begin{lemma}[Finite Symbolic State Space]
The symbolic state space is finite (see \hyperref[tab:fsm]{Table~\ref*{tab:fsm}}):
\[
S := \{0, \dots, 9\} \times \{0, 1\} \times \{0, 1, 2\}, \quad \text{so } |S| = 60
\]
\end{lemma}

\begin{remark}[Effective State Space is Reduced]
While the total symbolic state space contains $|S| = 60$ states,
only $\sim 20$ of them occur repeatedly during the symbolic evolution.
All states with $c > 0$ are restricted to the first symbolic transition,
and vanish immediately thereafter.
Thus, the effective dynamical system unfolds almost entirely within the set $\{(r, p, 0)\}$.
\end{remark}

\begin{figure}[H]
    \centering
    \includegraphics[width=0.85\textwidth]{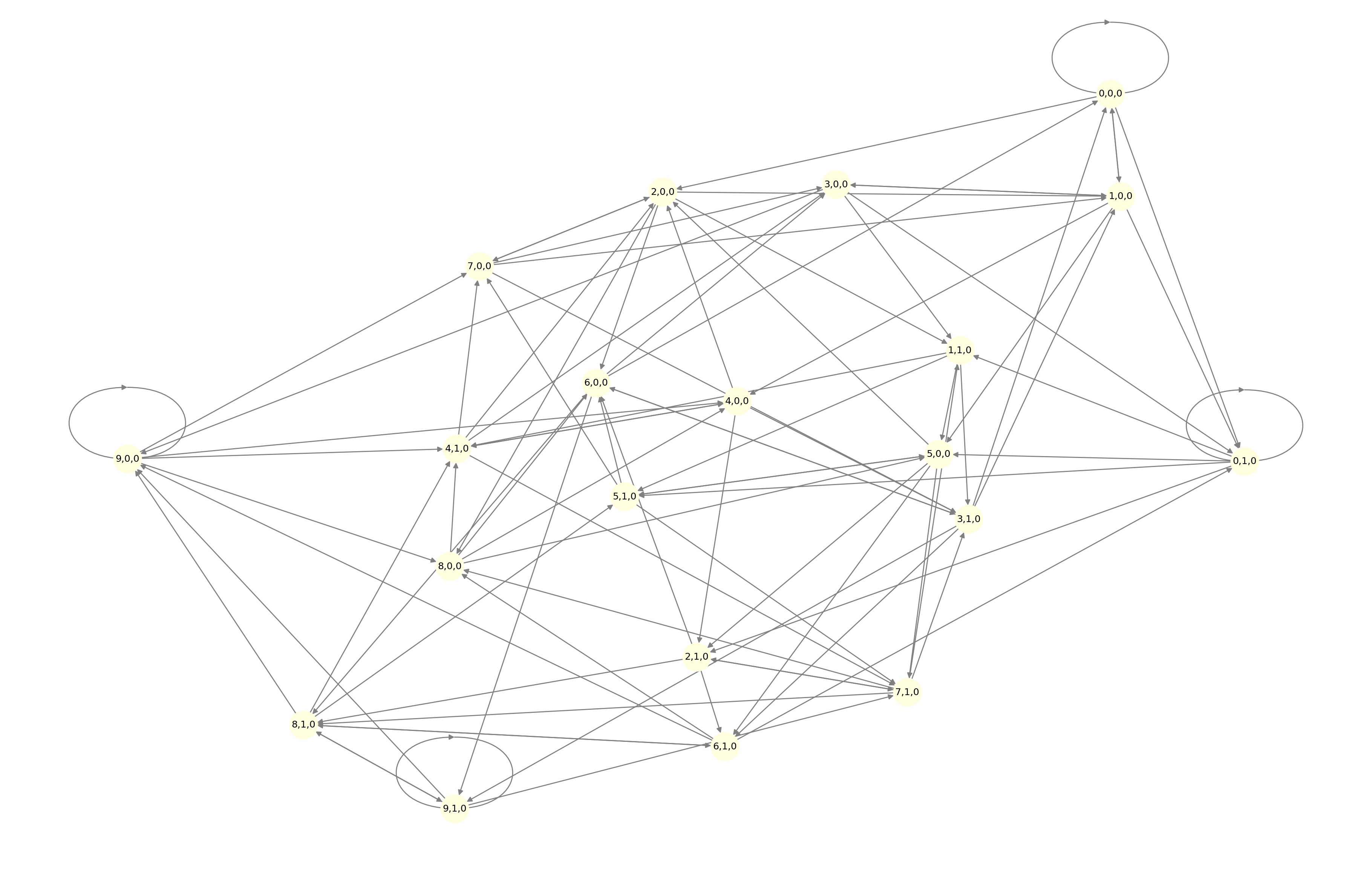}
    \caption{Transition graph of the symbolic FSM restricted to the 20 active states.  
    This subgraph contains the subset of symbolic triples \((r, p, c)\) that appear during actual Collatz iterations.  
    The transitions encode local digitwise updates under the symbolic emulation of \( T(n) \).}
    \label{fig:fsm20}
\end{figure}

\begin{theorem}[Terminal Symbolic Behavior and Global Absorption]
\label{thm:cycle}
Let $S := \{0,\dots,9\} \times \{0,1\} \times \{0,1,2\}$ denote the symbolic digit state space.

\begin{enumerate}
    \item Every symbolic state $(r, p, c) \in S \setminus \{(1,0,0), (2,0,0), (4,0,0)\}$ transitions deterministically to the null state $(0,0,0)$ under repeated application of the update rule $\delta$. 
    \item The three states $(1,0,0)$, $(2,0,0)$, and $(4,0,0)$ form a unique terminal cycle:
    \[
    (1,0,0) \to (4,0,0) \to (2,0,0) \to (1,0,0)
    \]
    However, this cycle is only stable when it occurs at the least significant state position $s_0$, and all higher positions $s_i$ with $i > 0$ have already reached the null state $(0,0,0)$.

    \item If any higher digit $s_k^{(t)} \neq (0,0,0)$ for some $k > 0$, then the terminal cycle is disrupted, and $s_0$ will continue evolving until the full collapse occurs.
\end{enumerate}
Thus, the symbolic Collatz system converges globally: every trajectory collapses into the null state in all but the least significant digit, which then cycles within the terminal loop.
\end{theorem}

\subsection{Local Symbolic Loops and Structural Collapse} \label{sec:local-loops}

While analyzing the symbolic FSM, we observe recurring local loops involving certain digit states. Examples include:
\begin{align*}
    &(6, 0, 0) \rightarrow (3, 0, 0) \rightarrow (0, 1, 0) \rightarrow (5, 1, 0) \rightarrow (6, 0, 0), \\
    &(8, 1, 0) \leftrightarrow (9, 1, 0), \\
    &(4, 1, 0) \rightarrow (7, 0, 0) \rightarrow (2, 0, 0) \rightarrow (1, 1, 0) \rightarrow (4, 1, 0).
\end{align*}

These finite cycles represent local oscillations induced by alternating parity patterns in the least significant state \( s_0 \). However, they do not persist globally. Any activation of a higher digit (e.g., through carry overflow during odd steps) disrupts the cycle and leads to symbolic absorption into the null state \( (0, 0, 0) \).

\begin{remark}[Local Cycles vs. Global Collapse]
Finite loops in the symbolic automaton do not contradict global convergence. They are:
\begin{itemize}
    \item confined to the least significant state \( s_0 \),
    \item dependent on specific combinations of parity and carry values,
    \item always disrupted once higher digits become active.
\end{itemize}
Hence, these loops are structurally harmless and inherently short-lived.
\end{remark}

\begin{corollary}[No Infinite Symbolic Loops Beyond the Terminal Cycle]
\label{cor:no-infinite-loops}
The only globally persistent symbolic cycle is the terminal 3-cycle
\[
(1, 0, 0) \rightarrow (4, 0, 0) \rightarrow (2, 0, 0) \rightarrow (1, 0, 0),
\]
which survives exclusively at state position \( s_0 \), after all higher digits have collapsed to \( (0, 0, 0) \).  
All other symbolic loops are local and inevitably collapse.
\end{corollary}

\paragraph{Loop-Driven Expansion and Role of \(s_0 \).}

Although these symbolic loops are transient, they play a key role in triggering symbolic expansion. In particular:

\begin{itemize}
    \item Loops in \( s_0 \) often alternate between odd and even states, producing a high frequency of symbolic transitions.
    \item This sustained alternation does \emph{not} require activity in higher digits—parity alone controls the switching.
    \item Over time, the rapid alternation can lead to carry buildup during odd steps (e.g., via \( 3r + c \)), which eventually activates \( s_1 \).
    \item Once higher digits are active, the loop breaks, and the symbolic system expands temporarily before collapsing.
\end{itemize}

\begin{remark}[Local timing and control element]
The least significant state \( s_0 \) acts as a symbolic \emph{clock generator}:  
its parity alternation determines the sequence of updates (odd vs. even), regardless of higher digit states.  
This rhythmic switching can sustain local loops, initiate expansion phases, and regulates the timing of symbolic contraction.
\end{remark}

\subsection{Stability of the Null State}

\begin{lemma}[Null state is absorbing]\label{lem:null-absorbing}
Let $s_i^{(t)} = (0,0,0)$ for some $i > 0$, and assume that the incoming carry is zero,
i.e., $c_i^{(t)} = 0$. Then $s_i^{(t+1)} = (0,0,0)$.
\begin{proof}
With $r = 0$ and $c = 0$ we have both $3r + c = 0$ (odd rule) and
$\lfloor(r - c)/2\rfloor = 0$ (even rule). Hence the symbolic state remains $(0,0,0)$.
\end{proof}
\end{lemma}

\begin{lemma}[Digit absorption]\label{lem:digit-absorption}
For every digit position $i > 0$, there exists $T_i \in \mathbb{N}$ such that
$s_i^{(t)} = (0,0,0)$ for all $t \ge T_i$.
\begin{proof}
We proceed by induction on $i$.

\emph{Base case ($i=1$):}
A newly created most significant digit can only arise through carry propagation
during an odd step. By the update rules defined in Section~\ref{sec:transition}, 
the only symbolic states that may appear in position $i=1$ are $(1,0,0)$ or $(2,0,0)$, 
resulting from a maximal carry of $2$. These states are part of the terminal cycle 
and collapse to $(0,0,0)$ in finitely many steps.

\emph{Inductive step:}
Assume that $s_i^{(t)} = (0,0,0)$ and $c_i^{(t)} = 0$ for all $t \ge T_i$.
Then, by Lemma~\ref{lem:null-absorbing}, the next higher digit satisfies
$s_{i+1}^{(t)} = (0,0,0)$ for all $t \ge T_i + 1$.

\emph{Conclusion:}
All digit positions $i > 0$ are eventually absorbed into the null state $(0,0,0)$,
and only the least significant digit remains dynamically active.
\end{proof}
\end{lemma}

\begin{theorem}[Finite symbolic support]\label{thm:finite-support}
For every symbolic trajectory $(s^{(t)})_{t \in \mathbb{N}}$, there exists a time $T \in \mathbb{N}$ such that
\[
s^{(t)} = [s_0^{(t)}, (0,0,0), (0,0,0), \dots] \quad \text{for all } t \ge T.
\]
Moreover, all transient higher-digit states are limited to the forms $(1,0,0)$ or $(2,0,0)$.
\end{theorem}

\begin{proof}
Immediate from Lemma~\ref{lem:digit-absorption}.
Combined with Theorem~\ref{thm:cycle}, this implies that after time $T$,
the complete symbolic configuration reduces to the 3-state cycle
\[
(1,0,0) \to (4,0,0) \to (2,0,0) \to (1,0,0)
\]
at the least significant state $s_0$.
\end{proof}

\subsection{Digit Emergence via Carry Propagation}

New digit positions in the decimal representation of $T(n)$ can only emerge through carry overflow at the most significant digit, typically during odd transitions of the form $T(n) = 3n + 1$. These overflows propagate leftward and may increment a previously zero digit, effectively increasing the overall length of the number.

However, this form of digit growth is structurally limited. Carry values depend on the alternating parity structure of the Collatz sequence, which itself is governed by the presence of trailing $1$-bits in the binary representation of the input. As these bits are progressively shifted out through repeated divisions, the system undergoes fewer odd transitions, and carry production diminishes accordingly.

\begin{remark}[Carry-Driven Digit Growth Is Self-Limiting]
Digit growth relies entirely on nonzero carry values originating from lower positions. Since such carries are only generated during odd steps—and these in turn depend on trailing $1$s in the binary expansion—they cannot persist indefinitely. As a result, the number of digits can increase only temporarily and eventually returns to a stable size.
\end{remark}

\subsection{Structural Constraint on the Peak Position}

We define the following key parameters of a symbolic Collatz orbit:

\begin{itemize}
  \item $t_{\text{end}}$: total orbit length until convergence to 1,
  \item $L_{\text{peak}}$: maximum number of decimal digits observed during the orbit,
  \item $t_{\text{peak}}$: first step at which this maximum digit length is reached.
\end{itemize}

From the definition of the effective symbolic decay time, we obtain:
\[
T_{\text{eff}} := \frac{t_{\text{end}} - t_{\text{peak}}}{L_{\text{peak}} - L_{\text{end}}}
\]
In the standard case, where $L_{\text{end}} = 1$, this simplifies to:
\[
T_{\text{eff}} := \frac{t_{\text{end}} - t_{\text{peak}}}{L_{\text{peak}} - 1}
\]

Rearranging this inequality gives a direct bound on $t_{\text{peak}}$:
\[
t_{\text{peak}} = t_{\text{end}} - T_{\text{eff}} \cdot (L_{\text{peak}} - 1)
\]

Using the trivial lower bound $T_{\text{eff}} \geq 1$ yields:
\[
t_{\text{peak}} \leq t_{\text{end}} - (L_{\text{peak}} - 1)
\]

And from the empirical estimate $T_{\text{eff}} \geq \frac{t_{\text{end}}}{L_{\text{peak}} - 1} - 2$,  
we finally derive the clean structural constraint:
\[
t_{\text{peak}} \leq 2(L_{\text{peak}} - 1)
\]

\paragraph{Interpretation.} The peak digit length is always reached sufficiently early in the orbit,  
specifically no later than twice the adjusted peak size.  
This constraint prevents excessive growth from occurring arbitrarily late in the trajectory  
and reflects an implicit structure in the digitwise evolution.

\paragraph{Implication.} The inequality reinforces the view that Collatz growth is globally constrained:  
even when large numeric peaks occur, the system compensates with a proportionally long collapse phase.  
Combined with the deterministic symbolic collapse (Theorem~4.8), this supports the hypothesis  
that infinite symbolic expansion is impossible.

\subsection{Symbolic Termination via Ranking Function}

\begin{definition}[Cycle Distance Ranking]
Let \( \mathcal{C} \subset S \) be the 3-cycle defined in Theorem~\ref{thm:cycle}.  
Define the symbolic ranking function \( \varrho : S \to \mathbb{N} \) by
\[
\varrho(s) := \min \{k \in \mathbb{N}_0 \mid \delta^{(k)}(s) \in \mathcal{C}\}
\]
\end{definition}

\begin{lemma}[Monotonic Descent]
Every symbolic transition $\delta(s_i) = s_{i+1}$ satisfies:
\[
\varrho(s_{i+1}) \leq \varrho(s_i),
\]
with strict inequality unless $s_i \in \mathcal{C}$
\end{lemma}

\begin{theorem}[Symbolic Termination]
All symbolic trajectories $(s_0, s_1, \dots)$ generated by iterated application of $\delta$ reach the terminal cycle $\mathcal{C}$ in finitely many steps.
\end{theorem}

\begin{proof}
The function $\varrho$ is finite-valued and strictly decreases outside $\mathcal{C}$. Thus, every sequence must eventually enter the cycle, and no infinite descent is possible.
\end{proof}

\begin{corollary}[Convergence of Collatz Iterates]
\label{cor:symbolic-termination}
For every $n \in \mathbb{N}^+$, the symbolic trajectory of $n$ under $\delta$ terminates in the 3-cycle $\mathcal{C}$, and hence the classical Collatz sequence enters the loop $4 \to 2 \to 1$.
\end{corollary}

\section{Binary Structural Convergence}
\label{sec:binary-convergence}

This section presents a binary-level framework for analyzing Collatz convergence.  
We decompose trajectories into symbolic growth and decay phases based on the trailing bits in the binary encoding.  
By isolating multiplicative and divisive behavior through tail structure, we derive both empirical and deterministic measures of convergence pressure.

\paragraph{Symbolic Block Formulation.}

We define the Collatz transformation in symbolic block form as follows:

\[
T(n) =
\begin{cases}
T_3^{\mathrm{to}(n)}(n), & \text{(multiplicative phase)} \\
T_2^{\mathrm{tz}(n')}(n'), & \text{(divisive phase), where } n' = T_3^{\mathrm{to}(n)}(n)
\end{cases}
\]

with:

\begin{align*}
T_3(n) &:= \frac{3n + 1}{2} \quad \text{(single odd step followed by division)}, \\
T_2(n) &:= \frac{n}{2} \quad \text{(pure halving step)}, \\
\mathrm{to}(n) &:= \text{number of trailing ones in binary of } n, \\
\mathrm{tz}(n) &:= \text{number of trailing zeros in binary of } n.
\end{align*}

Hence, each Collatz orbit proceeds in blocks of the form:

\[
n \to T_3^{\mathrm{to}(n)}(n) =: n' \to T_2^{\mathrm{tz}(n')}(n') =: n''
\]

This decomposition emphasizes the symbolic alternation between expansion and contraction, structured by the binary tail.

\subsection{Parity Alternation and Symbolic Lifetime}

Symbolic digit expansion under repeated Collatz steps requires alternating parity
to activate both the multiplication and division transitions in a controlled fashion.
This alternation arises from trailing $1$s in the binary structure of the input:
\[
\texttt{...1111}_2\quad\Rightarrow\quad\texttt{odd} \to \texttt{even} \to \texttt{odd} \to \cdots
\]

Each even step shifts the binary suffix rightward, replacing $1$s with $0$s.
As the tail of the binary representation flattens, parity changes become less frequent,
and the system eventually enters a stable shrinking phase dominated by divisions.

This behavior is illustrated by inputs such as \(n = 31 = \texttt{0b11111}\), which trigger long alternating sequences:
\[
31 \to 94 \to 47 \to 142 \to 71 \to \dots
\]
Each odd value is followed by an even one, until the bit suffix shifts:
\[
\texttt{...111} \to \texttt{...1000} \to \texttt{...01000} \to \dots
\]

\begin{remark}[Size Alone Does Not Imply Slowness]
It is a common misconception that large input values necessarily require the most steps to reach the terminal cycle. While large numbers can produce long orbits, symbolic slowdown is more strongly correlated with the \emph{binary structure} of the input.

This was already observed by Lagarias~\cite{lagarias1985generalizations}, who showed that numbers of the form $2^t - 1$ (e.g., $n = 31 = \texttt{0b11111}$) undergo $t$ consecutive odd steps before any division occurs. The trailing 1-bits in their binary encoding generate prolonged alternations, keeping symbolic length high despite numeric decrease.

As a result, such inputs exhibit the longest symbolic lifetimes per digit, making them the true worst-case cases for symbolic contraction---not necessarily the largest inputs.
\end{remark}

\begin{remark}[Binary Structure Governs Symbolic Lifetime]
Even extremely large numbers, such as \(10^{1000}-1\), often exhibit shorter symbolic durations than structured binary inputs like \(2^t - 1\).  
Thus, symbolic peak persistence depends more on bit entropy and parity alternation than numeric magnitude.
\end{remark}

\subsection{Binary Tail Energy Model and Structural Regulation}

In addition to the symbolic digitwise automaton, we can view the convergence dynamics through the lens of binary structure, particularly the number of trailing zeros in the binary expansion of $n$ at each step.

\begin{definition}[Binary Trailing Zero Count]
Let $n \in \mathbb{N}^+$. Define $\texttt{tz}(n)$ as the number of trailing $0$-bits in the binary representation of $n$:
\[
\texttt{tz}(n) := \max \{ k \in \mathbb{N} : 2^k \mid n \}
\]
\end{definition}

This count reflects how many consecutive even divisions $T(n)$ can undergo before encountering another odd value.

\paragraph{Interpretation.}
\begin{itemize}
    \item When $\texttt{tz}(n) = 0$, the Collatz update $T(n) = 3n + 1$ applies, and produces an even result: $\texttt{tz}(T(n)) \geq 1$.
    \item Each successive division by 2 reduces $\texttt{tz}$ by 1, until it reaches zero again.
\end{itemize}

This induces a natural alternation:
\[
\texttt{tz}=0 \Rightarrow \text{odd step (3n+1)} \quad \Rightarrow \quad \text{energy buildup}
\]
\[
\texttt{tz}>0 \Rightarrow \text{even steps} \quad \Rightarrow \quad \text{energy dissipation}
\]

\paragraph{Energy Interpretation.} 
We interpret $\texttt{tz}$ as a symbolic "energy buffer":
\begin{itemize}
    \item $3n+1$ increases symbolic digit length and triggers carry propagation (cf. Section~\ref{sec:transition}).
    \item Even steps reduce both the binary length and symbolic digit complexity.
    \item Thus, the system alternates between **growth** and **contraction**, but is globally biased toward contraction.
\end{itemize}

\paragraph{Structural Implication.}
This view supports the convergence theorem (Theorem~\ref{thm:finite-support}) from a binary angle: the symbolic system cannot grow arbitrarily, because each growth pulse from $3n+1$ is structurally paired with a trailing cascade of divisions that lower digit length and carry levels.

\begin{remark}[Terminal Energy Depletion]
Inputs like $n = 2^t - 1$ maximize symbolic alternation: their binary suffix consists of all $1$s. This forces $\texttt{tz} = 0$ for $t$ steps, leading to maximal symbolic expansion. But even these inputs eventually generate trailing zeros, and the system enters a decay phase.
\end{remark}

\begin{remark}[Bit-Length Growth Under \(3n + 1\) is Structurally Bounded]
The step \( T(n) = 3n + 1 \) increases the binary length of \( n \) by at most two bits. Specifically:

\[
\text{If } n \in [2^k, 2^{k+1}), \text{ then } 3n + 1 < 3 \cdot 2^{k+1} = 2^{k+1} \cdot 3
\Rightarrow \log_2(3n + 1) < k + 2
\]

Hence:
\[
\Delta_{\text{bits}}(T(n)) \leq 2
\]

This aligns with the decimal symbolic model, where a carry of at most 2 propagates during multiplication.

Since each even step (division by 2) reduces the bit length by exactly one, the growth induced by any odd step is compensated after at most two divisions.

This reinforces the interpretation that \texttt{tz} (the number of trailing zeros in binary) encodes a form of symbolic "energy" that the system discharges through repeated halving.
\end{remark}

\begin{table}[H]
\centering
\begin{tabular}{|c|r|c|c|c|}
\hline
\textbf{Step} & \textbf{Value} & \textbf{Binary} & \textbf{Bits} & \textbf{Trailing Zeros (tz)} \\
\hline
0  & 23   & \texttt{10111}      & 5  & 0 \\
1  & 70   & \texttt{1000110}    & 7  & 1 \\
2  & 35   & \texttt{100011}     & 6  & 0 \\
3  & 106  & \texttt{1101010}    & 7  & 1 \\
4  & 53   & \texttt{110101}     & 6  & 0 \\
5  & 160  & \texttt{10100000}   & 8  & 5 \\
6  & 80   & \texttt{1010000}    & 7  & 4 \\
7  & 40   & \texttt{101000}     & 6  & 3 \\
8  & 20   & \texttt{10100}      & 5  & 2 \\
9  & 10   & \texttt{1010}       & 4  & 1 \\
10 & 5    & \texttt{101}        & 3  & 0 \\
11 & 16   & \texttt{10000}      & 5  & 4 \\
12 & 8    & \texttt{1000}       & 4  & 3 \\
13 & 4    & \texttt{100}        & 3  & 2 \\
14 & 2    & \texttt{10}         & 2  & 1 \\
15 & 1    & \texttt{1}          & 1  & 0 \\
\hline
\end{tabular}
\caption{Bit-level Collatz iteration for $n = 23$}
\label{tab:collatz-23}
\end{table}

\begin{remark}[Bit-Level and Carry Coherence]
It is no coincidence that both the binary bit-length and the symbolic carry values are bounded by 2 during odd steps (i.e., $3n + 1$).  
This reflects a deep structural limit of the transformation:

\begin{itemize}
    \item \textbf{Binary growth:} Multiplication by 3 adds at most $\log_2(3) \approx 1.58$ bits. Thus,
    \[
    \texttt{bitlen}(3n + 1) \leq \texttt{bitlen}(n) + 2
    \]

    \item \textbf{Symbolic carry:} In the digitwise FSM, local updates $3r + c$ never exceed 29, so the propagated carry value satisfies
    \[
    c' = \left\lfloor \frac{3r + c}{10} \right\rfloor \leq 2 \quad \text{for all } r \in \{0, \dots, 9\},\, c \in \{0,1,2\}
    \]
\end{itemize}

This dual constraint---symbolic and binary---shows that the FSM encoding faithfully captures the core growth limits of the Collatz function.  
It is therefore not a numerical coincidence, but an algebraic invariant of the dynamics.
\end{remark}

\subsection{Bit-Length Bounds for Classical Collatz Growth}

\begin{theorem}[Bit-Length Growth Bound for Odd Steps]
\label{thm:bit-growth}
Let \( n \in \mathbb{N}^+ \), and let \( \beta(n) := \lfloor \log_2(n) \rfloor + 1 \) denote the bit-length of \( n \).  
Then the Collatz update step \( T(n) = 3n + 1 \), applied to odd \( n \), satisfies the bound
\[
\beta(T(n)) \leq \beta(n) + 2.
\]
\end{theorem}

\begin{proof}
For all \( n \in \mathbb{N}^+ \), we have:
\[
3n + 1 < 4n.
\]
Taking logarithms:
\[
\log_2(3n + 1) < \log_2(4n) = \log_2(n) + 2.
\]
Thus, using the bit-length function \( \beta(n) = \lfloor \log_2(n) \rfloor + 1 \), we get:
\[
\beta(3n + 1) = \lfloor \log_2(3n + 1) \rfloor + 1 \leq \lfloor \log_2(n) + 2 \rfloor + 1.
\]

Now two cases:

\textbf{Case 1:} \( \log_2(n) \in \mathbb{Z} \).  
Then \( \lfloor \log_2(n) + 2 \rfloor = \log_2(n) + 2 \), so
\[
\beta(3n + 1) \leq \log_2(n) + 2 + 1 = \beta(n) + 2.
\]

\textbf{Case 2:} \( \log_2(n) \notin \mathbb{Z} \).  
Then \( \lfloor \log_2(n) + 2 \rfloor = \lfloor \log_2(n) \rfloor + 2 \), and again:
\[
\beta(3n + 1) \leq \lfloor \log_2(n) \rfloor + 2 + 1 = \beta(n) + 2.
\]

In both cases, the bit-length increases by at most 2.
\end{proof}

\subsection{Two-Phase Decomposition}

Starting from an odd number \( n \), the Collatz sequence can be decomposed into two symbolic phases:

\paragraph{Phase 1 – Multiplicative Phase (\( T_3 \)-block):}

\[
T_3(T_3(\dots T_3(n))) \quad \text{(applied } t_o(n) \text{ times)},
\]

This corresponds to a block of \( t_o(n) \) successive applications of the transformation \( T(n) = \frac{3n + 1}{2} \). Each such step tends to increase the bit-length by approximately \( \log_2\left(\frac{3}{2}\right) \).

\paragraph{Phase 2 – Divisive Phase (\( T_2 \)-block):}

Let \( n' \) be the result of the \( T_3 \)-block. Then:

\[
T_2(T_2(\dots T_2(n'))) \quad \text{(applied } t_z(n') \text{ times)},
\]

This phase corresponds to successive halving steps (i.e., \( T(n) = \frac{n}{2} \)), each reducing the bit-length by exactly one.

\subsection{Interpretation of Terms}

\begin{itemize}
  \item \( t_o(n) \): the number of trailing ones in the binary representation of \( n \)
  \item \( n' \): the result of applying the transformation \( T(n) = \frac{3n + 1}{2} \) exactly \( t_o(n) \) times
  \item \( t_z(n') \): the number of trailing zeros in the binary representation of \( n' \)
\end{itemize}

\subsection{Symbolic Drift Function}
\label{sec:symbolic-drift}

We define the symbolic drift function as
\[
w(n) := \mathrm{to}(n) \cdot \log_2\left(\tfrac{3}{2}\right) - \mathrm{tz}(n'),
\]
where \( \mathrm{to}(n) \) is the number of successive symbolic \( T_3 \)-steps starting from an odd number \( n \), and \( \mathrm{tz}(n') \) is the number of trailing zeros in binary of the result \( n' := T_3^{\mathrm{to}(n)}(n) \). This function measures the net symbolic bit-length change of a Collatz block and reflects the balance between expansion and contraction.

\paragraph{Worst-case instance.}
For inputs of the form \( n = 2^k - 1 \), we have \( \mathrm{to}(n) = k \) and typically \( \mathrm{tz}(n') \geq k \), yielding
\[
w(n) \leq k \cdot \left( \log_2\left(\tfrac{3}{2} \right) - 1 \right) = k \cdot \log_2\left(\tfrac{3}{4}\right) < 0,
\]
matching the symbolic contraction constant \( \log_2(3/4) \approx -0.415 \) noted by Tao.

\paragraph{Empirical average.}
Simulations over all odd \( n \leq 10^6 \) yield:
\[
\mathbb{E}[w(n)] \approx -0.83008,
\]
confirming that symbolic contraction dominates on average. This value appears consistently across large input ranges and supports a universal convergence bias.

\begin{figure}[H]
    \centering
    \includegraphics[width=\textwidth]{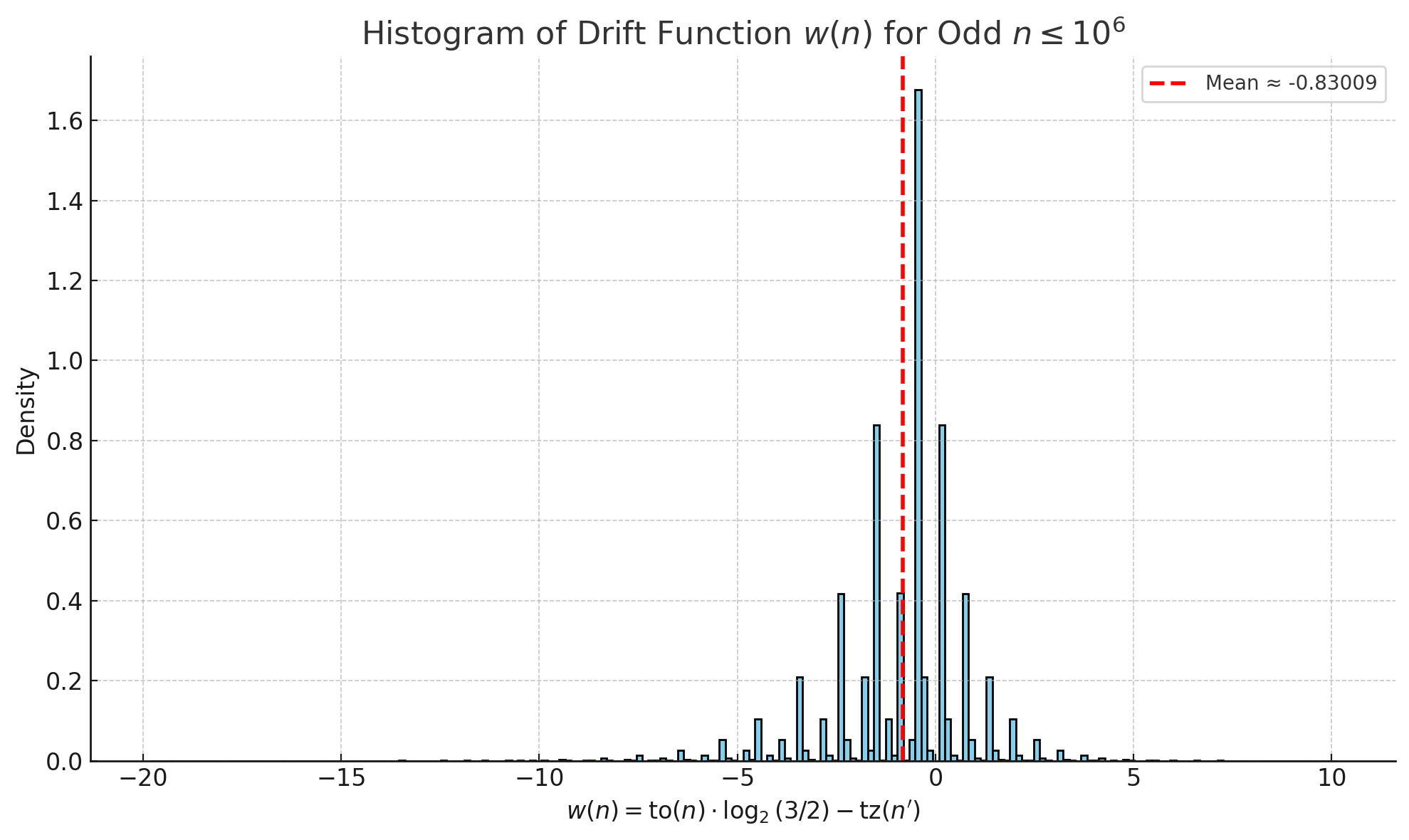}
    \caption{Distribution of symbolic drift $w(n)$ for all odd $n \leq 10^6$. The mean value (dashed red) confirms a consistent structural contraction.}
    \label{fig:drift-distribution}
\end{figure}

\begin{remark}[Gaussian Profile and Central Limit Behavior of Symbolic Drift]
When computed over all odd inputs \( n \leq 10^6 \), the symbolic drift function \( w(n) \) exhibits a nearly Gaussian distribution centered around \( \mathbb{E}[w(n)] \approx -0.83 \). The empirical histogram shows a symmetric bell shape with skewness \( \approx -1.09 \) and excess kurtosis \( \approx 4.03 \), indicating mild left-skew and slightly heavier tails compared to the standard normal curve.

This behavior is particularly remarkable given that \( w(n) \) is defined through purely deterministic symbolic transformations involving \( \mathrm{to}(n) \) and \( \mathrm{tz}(n') \). Despite their discrete origin, these terms aggregate statistically as if governed by the central limit theorem: they act like weakly dependent components whose combined effect produces a continuous, tightly concentrated distribution.

The resulting normality highlights an emergent regularity within symbolic Collatz dynamics -- a convergence structure that is not only average-negative but also statistically regulated. It supports the interpretation of \( w(n) \) as a structural invariant rather than a mere empirical artifact.
\end{remark}

\begin{figure}[H]
    \centering
    \includegraphics[width=0.75\textwidth]{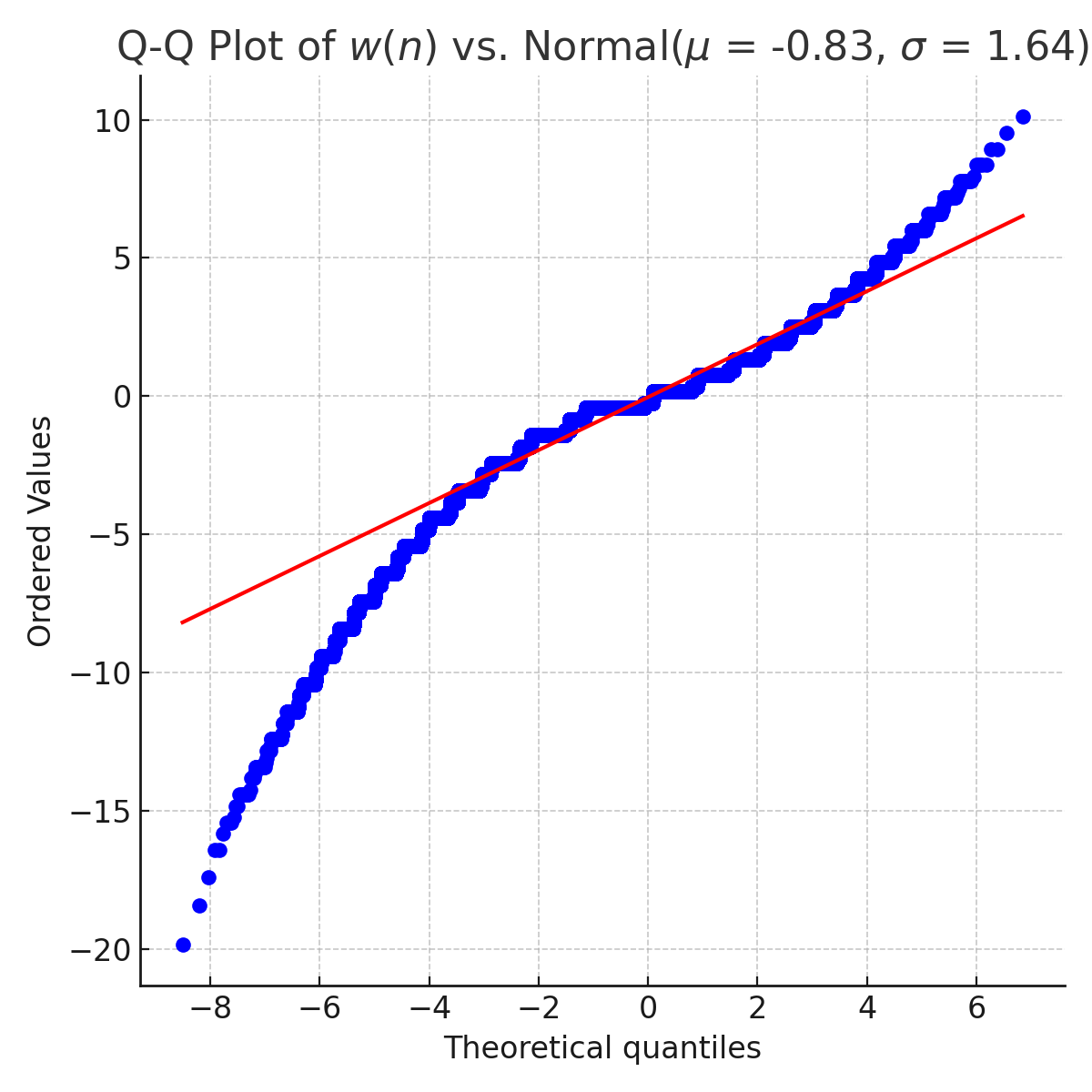}
    \caption{Q–Q plot comparing the empirical distribution of the symbolic drift function \( w(n) \) (for odd \( n \leq 10^6 \)) against a normal distribution with \( \mu = -0.83 \), \( \sigma \approx 1.638 \). The close alignment confirms near-Gaussian behavior.}
    \label{fig:qqplot}
\end{figure}

\begin{remark}[Q–Q-Plot and Empirical Normality of the Drift]
To empirically verify the near-normality of the drift distribution \( w(n) \), we compared it to a Gaussian distribution \( \mathcal{N}(\mu = -0.83, \sigma \approx 1.638) \) using a quantile–quantile (Q–Q) plot. The empirical quantiles align closely with the theoretical normal quantiles, particularly in the central region of the distribution.

Minor deviations near the tails confirm the slightly leptokurtic profile observed earlier (kurtosis \( \approx 4.03 \)) and a mild left skew (skewness \( \approx -1.09 \)). Nevertheless, the overall fit is strong, confirming that the symbolic drift behaves statistically like a Gaussian-distributed variable.

This alignment reinforces the interpretation of \( w(n) \) as an emergent statistical invariant: a deterministic but stochastically behaving quantity, supporting the idea that the symbolic Collatz dynamics simulate a central limit phenomenon without requiring randomness.
\end{remark}

\begin{remark}[Emergent Structure Behind Probabilistic Success]
The near-Gaussian behavior of the symbolic drift function \( w(n) \), despite its fully deterministic definition, suggests that the statistical patterns observed in previous stochastic or heuristic approaches may not arise from inherent randomness — but rather from the structured composition of symbolic Collatz dynamics itself.

This perspective reframes the effectiveness of probabilistic models (e.g., via random walks or ergodic arguments) not as evidence for randomness in the system, but as a reflection of underlying regularities that mimic stochastic behavior. In this light, previous results may have succeeded not by modeling true uncertainty, but by coincidentally aligning with emergent statistical invariants of a deterministic process.
\end{remark}

\paragraph{Structural explanation.}
This contraction is not accidental but structurally determined by binary properties of the input. Specifically:

\begin{enumerate}
  \item \textbf{Trailing-one expectation.}  
  For odd \( n \), the number of trailing binary ones satisfies:
  \[
  \mathbb{E}[\mathrm{to}(n)] = \sum_{k=1}^\infty k \cdot \left(\frac{1}{2}\right)^k = 2
  \]

  \item \textbf{Decay behavior.}  
  Empirically, the expansion result \( n' \) satisfies:
  \[
  \mathbb{E}[\mathrm{tz}(n')] \approx 2.00
  \]
  Thus, both phases—symbolic expansion and decay—have symmetric average length.

  \item \textbf{Resulting contraction.}  
  Substituting into the drift formula gives:
  \[
  \mathbb{E}[w(n)] \approx 2 \cdot \log_2\left(\tfrac{3}{2} \right) - 2 = \log_2\left( \tfrac{9}{16} \right) \approx -0.83008
  \]
\end{enumerate}

\paragraph{Significance.}
This value represents a structural invariant of the symbolic Collatz system. Each block enforces net bit-length loss, making long-term expansion impossible.

\begin{remark}[Role of Division Steps]
In symbolic Collatz blocks, the halving steps \( T_2(n) = \frac{n}{2} \) play a strictly subtractive role. Each application removes one trailing binary zero without altering internal structure:
\[
\text{Bin}(T_2(n)) = \text{Bin}(n) \text{ with the final } 0 \text{ removed}.
\]
Hence, \( \mathrm{tz}(n') \) contributes a pure contraction term in the drift function, structurally balancing the additive impact of \( T_3 \)-steps.
\end{remark}

\subsection{Illustrative Example: Symbolic Block Decomposition for \texorpdfstring{$n = 15$}{n = 15}}

We demonstrate the symbolic block decomposition for the Collatz orbit starting at $n = 15$, using the alternation structure defined by trailing ones ($\texttt{to}$) and trailing zeros ($\texttt{tz}$).

\paragraph{Full Orbit (values only):}
\[
15 \to 46 \to 23 \to 70 \to 35 \to 106 \to 53 \to 160 \to 80 \to 40 \to 20 \to 10 \to 5 \to 16 \to 8 \to 4 \to 2 \to 1
\]

\paragraph{Block Decomposition (step-based):}
Each \textit{step} refers to a single transition (an arrow $a \to b$).  
The value $n'$ at the end of a $T_3$-block is defined as the **final value reached before switching to a divisive phase**.  
The number of divisions is then governed by $\texttt{tz}(n')$.

\begin{itemize}
  \item \textbf{Block 1 – \boldmath$T_3$ (8 steps):}  
  Triggered by $\texttt{to}(15) = 4$, which leads to 4 symbolic $T_3$ cycles  
  (each consisting of a \texttt{mul3} and a \texttt{div2}):
  \[
  15 \to 46 \to 23 \to 70 \to 35 \to 106 \to 53 \to 160 \to 80
  \]
  The final value is $n' = 80$, so the next block will contain exactly $\texttt{tz}(80) = 4$ halving steps.

  \item \textbf{Block 2 – \boldmath$T_2$ (4 steps):}
  \[
  80 \to 40 \to 20 \to 10 \to 5
  \]
  These are pure division steps, ending when $\texttt{to}(5) = 1$ triggers a new symbolic cycle.

  \item \textbf{Block 3 – \boldmath$T_3$ (2 steps):}
  \[
  5 \to 16 \to 8
  \]
  Based on $\texttt{to}(5) = 1$, we apply one $T_3$ cycle.

  \item \textbf{Block 4 – \boldmath$T_2$ (3 steps):}
  \[
  8 \to 4 \to 2 \to 1
  \]
  Even though $\texttt{tz}(8) = 3$, the orbit ends after exactly 3 division steps.
\end{itemize}

\paragraph{Summary:}
\[
\boxed{T_3\text{-Block (8 steps)}} \;\to\; 
\boxed{T_2\text{-Block (4 steps)}} \;\to\; 
\boxed{T_3\text{-Block (2 steps)}} \;\to\; 
\boxed{T_2\text{-Block (3 steps)}}
\]

\begin{remark}[Precise Phase Boundaries]
The transition from a $T_3$ to a $T_2$ block occurs \textbf{after the final result} of the last $T_3$ step has been reached.  
This means that $n'$—used to determine the next $\texttt{tz}$—is defined as the terminal value of the $T_3$-phase.  
In this case: $n' = 80$, and $\texttt{tz}(80) = 4$.
\end{remark}

\subsection{Symbolic Energy Function and Deterministic Convergence}

To move from empirical trends to a formal descent guarantee, we define a symbolic potential function that strictly decreases along Collatz trajectories. Let \( n_0 \) denote the initial value of a given trajectory. We define:
\[
f(n) := \frac{\log_2(n)}{\log_2(n_0)} + \mathrm{to}(n) - \mathrm{tz}(n),
\]
where \( \mathrm{to}(n) \) is the number of trailing ones in the binary representation of \( n \), and \( \mathrm{tz}(n) \) the number of trailing zeros.

This function captures both the relative magnitude of the current iterate and its symbolic structure. It contracts whenever a symbolic block—comprising \( \mathrm{to}(n) \) many \( T_3 \)-steps followed by \( \mathrm{tz}(n') \) divisions—has net negative drift.

\begin{figure}[H]
    \centering
    \includegraphics[width=\textwidth]{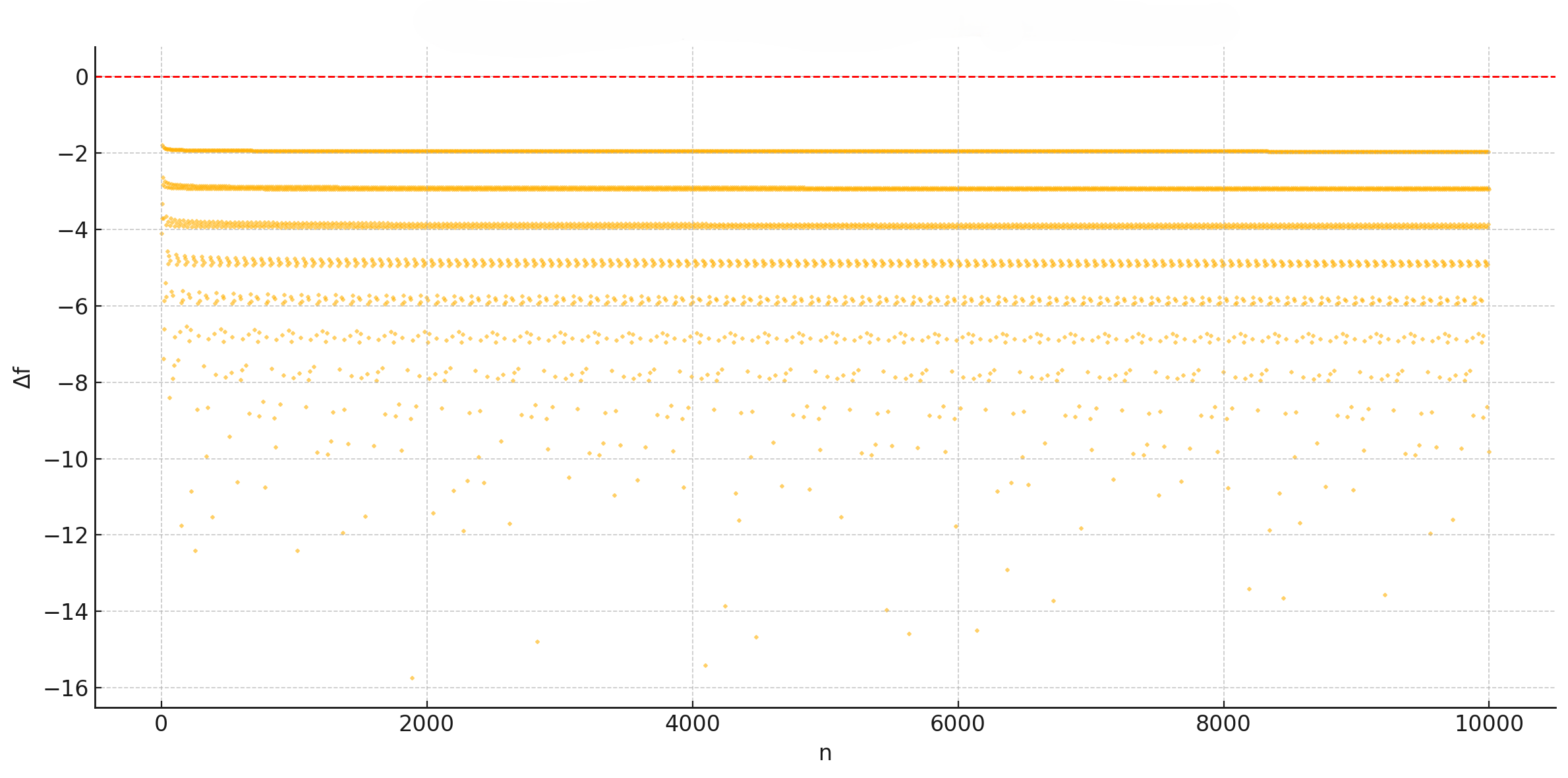}
    \caption{
    Change in the symbolic energy function over the range \( n \in [2, 10^4] \) for all inputs with \( \mathrm{to}(n) > 0 \).  
    Each point represents the difference \( \Delta f = f(T^{\mathrm{to}(n)}(n)) - f(n) \).  
    The consistent negative values confirm that symbolic Collatz blocks induce strict contraction of the ranking function.  
    This supports deterministic convergence along structured trajectories.
    }
    \label{fig:deltaf-plot}
\end{figure}
\paragraph{Descent Property.}
Let \( k := \mathrm{to}(n) \), \( m := \mathrm{tz}(n') \), and \( n' := T_3^k(n) \). Then:
\[
f(T^k(n)) - f(n) = \frac{\log_2(n') - \log_2(n)}{\log_2(n_0)} - (k + m)
\]
Using the estimate \( n' \leq (3/2)^k \cdot n \), we obtain:
\[
f(T^k(n)) - f(n) \leq \frac{k \cdot \log_2(3/2)}{\log_2(n_0)} - (k + m)
\]
Setting \( \log_2(n_0) = 1 \) for normalization, we recover the symbolic drift inequality:
\[
1 + \frac{\mathrm{tz}(n')}{\mathrm{to}(n)} > \log_2\left( \tfrac{3}{2} \right)
\]
which always holds, since \( \mathrm{to}(n) \geq 1 \) and \( \mathrm{tz}(n') \geq 0 \).

\begin{lemma}[Symbolic Energy Decrease]
For all odd \( n \) with \( \mathrm{to}(n) > 0 \), the symbolic energy function strictly decreases:
\[
f(T^{\mathrm{to}(n)}(n)) < f(n)
\]
\end{lemma}

\paragraph{Conclusion.}
The symbolic potential function \( f(n) \) provides a deterministic ranking that decreases blockwise and is bounded from below. Therefore, the symbolic Collatz system terminates without statistical assumptions. Each block enforces bit-level decay and guarantees symbolic contraction.

\begin{remark}[Absolute Drift Inequality]
The inequality
\[
1 + \frac{\mathrm{tz}(n')}{\mathrm{to}(n)} > \log_2\left(\tfrac{3}{2} \right)
\]
is not a statistical average but a symbolic invariant of every odd input.
\end{remark}

\begin{remark}[Historical Note on the Drift Constant]
The symbolic drift constant
\[
\log_2\left( \tfrac{3}{2} \right) - 1 = \log_2\left( \tfrac{3}{4} \right) \approx -0.415
\]
has appeared in earlier work:

\begin{itemize}
  \item Terras and Lagarias observed it in average-case and density models for Collatz iteration.
  \item Tao derived it in a probabilistic setting using stochastic dynamics.
\end{itemize}

However, the expression
\[
\log_2\left( \tfrac{3}{2} \right) - 1
\]
now admits a structural interpretation: one symbolic growth (\( \log_2(3/2) \)) followed by one contraction step (division by 2). This alternation explains the convergence behavior not through probability, but as a deterministic invariant encoded in the FSM.
\end{remark}

\subsection{Connection to Decimal FSM}

While the drift and energy functions are derived purely from the binary representation of \( n \), they align structurally with the decimal FSM model. 

\begin{itemize}
  \item Both frameworks reflect the same growth–decay duality: odd steps expand symbolic size, even steps contract it.
  \item The FSM enforces bounded carries (\( c \in \{0,1,2\} \)) during symbolic multiplication, which corresponds to a maximal bit-length increase of 2 under \( 3n+1 \).
  \item The convergence theorem for the FSM (Theorem~\ref{thm:finite-support}) is structurally mirrored by the symbolic drift inequality:
  \[
  1 + \frac{\mathrm{tz}(n')}{\mathrm{to}(n)} > \log_2\left( \tfrac{3}{2} \right),
  \]
  which guarantees descent in the binary energy function \( f(n) \).
\end{itemize}

This approach provides a new structural lens for analyzing recursive functions. Rather than relying on numerical growth bounds, it translates dynamics into a finite symbolic domain, making it amenable to mechanized reasoning. While the main construction operates over base-10 representations, the method itself is not inherently decimal: symbolic emulation has also been demonstrated for binary encodings using adapted FSM variants.

Two such binary FSM models are outlined in Appendix~\ref{app:binaryfsm}, showing how the symbolic transition idea extends beyond decimal digits. Whether this framework generalizes to broader function classes or arbitrary number bases remains an open question.

\section{Limitations and Formal Scope}

The symbolic Collatz framework developed in this paper replaces classical analytic or probabilistic strategies with a fully deterministic and constructive model. The digitwise automata introduced in Section~\ref{sec:transition} emulate the Collatz function \( T(n) \) through local updates on finite symbolic states.

\begin{itemize}
  \item All symbolic operations — digit parity, carry propagation, and the local update rules $\delta_{\text{odd}}$, $\delta_{\text{even}}$ — are elementary and primitive recursive. The entire system is, in principle, representable within first-order Peano Arithmetic (PA) or bounded arithmetic.

  \item The correctness of the digitwise emulation (Theorem~\ref{thm:encoding}) ensures that symbolic steps correspond exactly to Collatz updates. While this equivalence has not yet been formally verified in a proof assistant such as Coq or Lean, the symbolic definitions are fully explicit and mechanization-ready.

  \item The convergence proof (Theorem~\ref{thm:finite-support}) is conditional on the correctness of the encoding. Once this correspondence is formalized, the result becomes provable without reference to classical number theory or analytic bounds.
\end{itemize}

This symbolic approach offers a new structural lens for analyzing recursive functions. Rather than relying on global numerical estimates, it translates dynamics into a finite symbolic domain with provably bounded behavior, making it amenable to automated verification.

Although the main construction operates in decimal (base-10), the method is not inherently tied to this number system. In Appendix~\ref{app:binaryfsm}, we describe two alternative FSM variants built for binary (base-2) inputs. These models follow the same symbolic principles, but require modified local update rules due to the differences in base arithmetic:

\begin{itemize}
  \item Carry propagation behaves differently in binary and leads to a reduced symbolic state space.
  \item Parity-based transitions simplify, but the symbolic encoding must account for bitwise alignment and mod-2 carries.
  \item As in the decimal case, the global arithmetic behavior of \( T(n) \) is reproduced through strictly local symbolic transitions.
\end{itemize}

Finally, although our framework was developed for the classical Collatz function, it generalizes naturally to arbitrary affine maps \( T(n) = a n + b \). As shown in Section~\ref{sec:anb-generalization}, the same digitwise mechanism can emulate such functions via adjusted local rules. This opens the door to broader classes of recursively defined systems whose dynamics can be captured symbolically.

\section{Discussion and Future Directions}

The symbolic finite-state model introduced in this work offers a new structural framework for analyzing arithmetic dynamics through local, digitwise transformations. While originally constructed for the classical Collatz function, its underlying architecture—based on parity-aware carry propagation—is fundamentally extensible.

\paragraph{Generalization Potential.}
The digitwise mechanism extends naturally to affine maps of the form \( T(n) = a n + b \), as long as the coefficients are compatible with local carry arithmetic. Appendix~\ref{sec:anb-generalization} outlines a generalized transition scheme for such updates. This suggests the possibility of a broader symbolic calculus for integer recursions, enabling new finite-state models beyond the Collatz case.

\paragraph{Binary Extensions and Structural Insights.}
In addition to the decimal representation, binary FSM variants (Appendix~\ref{app:binaryfsm}) implement the same symbolic principles in base 2. While these require different transition rules, they retain the key idea: reconstructing global behavior through strictly local symbolic updates. When combined with the drift and ranking functions of Section~\ref{sec:binary-convergence}, this binary view provides a structural explanation for convergence, independent of numerical bounds.

\paragraph{Open Problems.}
The finite-state framework opens several directions for further exploration:
\begin{itemize}
  \item Can symbolic automata be constructed for other unsolved problems in arithmetic dynamics?
  \item How do symbolic trajectories behave under variations of \( a \) and \( b \) in affine maps \( T(n) = a n + b \)?
  \item Can symbolic convergence principles be extended to non-affine or piecewise-defined recursive functions?
\end{itemize}
While base changes (e.g., decimal to binary) alter local update rules, the overall convergence structure appears invariant—pointing to a deeper algebraic principle that merits further investigation.

\paragraph{Claim of Novelty.}
To the best of the author's knowledge, this is the first model that emulates Collatz dynamics via explicit symbolic digit triples \((r, p, c)\), deterministic local transitions, and bounded carry propagation within a finite-state system. Unlike heuristic or analytic approaches, the present framework offers a constructive reformulation of the problem, enabling formal convergence proofs and laying the groundwork for mechanized verification and broader generalization.

\section{Conclusion and Structural Insights}

We have introduced a finite-state symbolic automaton that emulates the Collatz map \( T(n) \) through strictly local digitwise transformations. Each decimal digit is represented by a symbolic triple \((r, p, c)\), encoding the digit value, the parity of the next digit, and a bounded carry. These local updates—total, deterministic, and primitive recursive—collectively reconstruct the full Collatz iteration without global arithmetic, enabling structural convergence proofs within a finite symbolic space.

\paragraph{Symbolic Convergence.}
The automaton enforces termination by collapsing all but the least significant symbolic digit into the null state \((0,0,0)\). Once only \( s_0 \) remains active, the system enters a unique 3-cycle \((4,0,0) \rightarrow (2,0,0) \rightarrow (1,0,0)\), structurally mimicking the classical loop \(4 \to 2 \to 1\). Convergence is guaranteed by bounded carry propagation and digit-level absorption, showing that the automaton supports only finitely many dynamically active states throughout any trajectory.

\paragraph{Binary Interpretation.}
To complement the symbolic model, we introduced a binary framework that analyzes convergence through the number of trailing ones and zeros in binary expansion. This yields a structural drift function and an energy ranking potential that contract deterministically across symbolic blocks. Odd steps \( (3n + 1)/2 \) increase symbolic complexity and bit length, while halving steps remove trailing zeros without altering internal structure. Crucially, both the binary bit-length and the symbolic carries are strictly bounded—odd steps increase size by at most two bits, and symbolic carries never exceed two. This confirms a fundamental dissipative property embedded in the dynamics.

\paragraph{Dual Perspective.}
Taken together, the symbolic and binary views reveal that Collatz convergence is not merely numerical but deeply structural: a self-limiting alternation between symbolic growth and digit-level decay. The apparent complexity of Collatz trajectories arises from a highly constrained interplay of parity, position, and bounded carry logic—far from chaotic, the process behaves like a regulated symbolic system with finite entropy.

\paragraph{Outlook.}
This framework opens a new pathway for analyzing recursive arithmetic functions via finite-state models. Future work may:
\begin{itemize}
  \item mechanize the symbolic encoding in proof assistants such as Coq or Lean,
  \item generalize the automaton to other affine maps \( T(n) = a n + b \) with symbolic carry logic,
  \item explore base-2 and modular variants with adapted symbolic states,
  \item or design related symbolic dynamics for other open problems in arithmetic iteration.
\end{itemize}

The symbolic Collatz FSM thus recasts a historically intractable conjecture into a verifiable, locally deterministic system. By shifting the focus from numerical iteration to symbolic structure, it reveals a new finite-state foundation for convergence—and opens the door to further exploration at the interface of number theory, logic, and automata.

\section*{Acknowledgements}

This work was developed independently. I would like to acknowledge the assistance of OpenAI’s language model (ChatGPT-4), which supported the writing process by helping with linguistic refinement, formal structuring, and programming tasks related to the generation of illustrative graphics. All mathematical reasoning, results, and conclusions are solely my own.

\vspace{1em}
\noindent
\textbf{License:} This work is licensed under a \href{https://creativecommons.org/licenses/by-nc-nd/4.0/}{CC BY-NC-ND 4.0} license.

\bibliographystyle{unsrt}
\bibliography{references}

\appendix
\section{Formal and Computational Details}

This appendix provides the formal underpinnings and supplementary computations that support the symbolic Collatz model developed in the main text. It contains:

\begin{itemize}
  \item a generalization of the digitwise transition mechanism to arbitrary affine functions \( T(n) = a n + b \),
  \item step-by-step symbolic evaluations for selected input values,
  \item a worked example of symbolic state evolution along a full Collatz trajectory, and
  \item the complete transition table of the finite-state machine governing symbolic digit dynamics.
\end{itemize}

All constructions are presented in a primitive recursive form and are fully compatible with the deterministic automaton architecture introduced in Section~\ref{sec:transition}.

\section{Symbolic Digitwise Model for General \boldmath{$an + b$}}

The digitwise update mechanism developed for the Collatz map \( T(n) = 3n + 1 \) generalizes naturally to arbitrary affine functions of the form
\[
T(n) = a \cdot n + b,
\]
with \( a, b \in \mathbb{N} \). This section describes how such functions can be emulated through purely local digit-level transformations. Two distinct strategies are available, depending on the size of the parameters \( a \) and \( b \):

\subsection{Digitwise Transition Rule}

Let \( n = \sum_{i=0}^L r_i \cdot 10^i \) be the decimal expansion of \( n \), and define an initial carry \( c_0 := 0 \). There are two strategies for computing \( T(n) = a \cdot n + b \) symbolically:

\paragraph{(1) Simplified Rule for Small \boldmath{$a,b$}:}
If \( a < 10 \) and \( b < 10 \), then the digitwise update can be performed directly during the multiplication phase:
\[
r_i' :=
\begin{cases}
(a \cdot r_i + b + c_i) \bmod 10 & \text{if } i = 0, \\
(a \cdot r_i + c_i) \bmod 10 & \text{if } i > 0,
\end{cases}
\qquad
c_{i+1} :=
\begin{cases}
\left\lfloor \dfrac{a \cdot r_i + b + c_i}{10} \right\rfloor & \text{if } i = 0, \\
\left\lfloor \dfrac{a \cdot r_i + c_i}{10} \right\rfloor & \text{if } i > 0.
\end{cases}
\]
This method is valid only for single-digit constants \( b \).

\paragraph{(2) General Two-Phase Method (Fully Correct):}
For arbitrary \( a, b \in \mathbb{N} \), the symbolic computation must be performed in two separate phases:
\begin{enumerate}
  \item \textbf{Multiplication:} Compute \( a \cdot n \) digitwise with carry propagation:
  \[
  r_i' := (a \cdot r_i + c_i) \bmod 10, \qquad c_{i+1} := \left\lfloor \frac{a \cdot r_i + c_i}{10} \right\rfloor
  \]
  \item \textbf{Addition:} Perform a standard digitwise addition of \( b \) to the result of the multiplication, including carry.
\end{enumerate}
This ensures correctness for all \( a, b \in \mathbb{N} \).

\subsection{Alternative Formulation}

The same transformation can be written as:
\[
T(n) = a n + b = \left( \sum_{i=0}^L a r_i \cdot 10^i \right) + b,
\]
which corresponds to multiplying each digit \( r_i \) by \( a \), applying carry propagation, and finally injecting the constant term \( b \) at the least significant position.

\begin{remark}[Full Addition of \boldmath{$b$}]
For general \( b \in \mathbb{N} \), the additive term should be treated as a separate value added to \( a \cdot n \) after digitwise multiplication and carry propagation. The simplified digitwise injection of \( b \) at the least significant digit (i.e., within \( r_0' \)) remains valid only if \( b < 10 \). For larger \( b \), a full digitwise addition must follow the multiplication phase.
\end{remark}

\subsection{Bounded Carry and Finite State Space}

The carry values remain bounded due to the fixed-digit range:
\[
c_{\max} \le \left\lfloor \frac{9a + b}{10} \right\rfloor,
\]
so the symbolic state space remains finite:
\[
S_{(a,b)} := \{0, \dots, 9\} \times \{0,1\} \times \{0, \dots, c_{\max}\}.
\]

\subsection{Interpretation}

This generalized rule preserves the core structure of the symbolic Collatz system: digitwise updates, carry propagation, and optional parity tagging. It allows the symbolic emulation of arbitrary affine maps of the form \( T(n) = a n + b \) via finite-state digit dynamics in base 10.

In this light, the classical Collatz rule \( T(n) = 3n + 1 \) becomes a special case of the more general symbolic framework. The model thus extends beyond Collatz-specific behavior and opens the door to structural analysis of other recursive systems governed by similar local rules.

\subsection{Generalized Symbolic Drift in \( T_{a,1} \)-Systems}

\begin{definition}[Symbolic Block Drift]
Let \( T_{a,1}(n) := \frac{a n + 1}{2^k} \), where \( k := \mathrm{tz}(a n + 1) \).  
We define the symbolic drift of a full Collatz block (growth–decay phase) as:
\[
w_{a,1}(n) := \mathrm{to}(n) \cdot \log_2\left( \frac{a}{2} \right) - \mathrm{tz}\left(T_3^{\mathrm{to}(n)}(n)\right)
\]
where:
\begin{itemize}
  \item \( \mathrm{to}(n) \): number of trailing ones in the binary expansion of \( n \),
  \item \( T_3(n) := \frac{a n + 1}{2} \): odd-case transformation,
  \item \( T_3^{\mathrm{to}(n)}(n) \): value after applying \( \mathrm{to}(n) \) many \( T_3 \)-steps,
  \item \( \mathrm{tz}(\cdot) \): number of trailing zeros (i.e., maximal halving steps following the expansion).
\end{itemize}
\end{definition}

\begin{remark}[Structural Contraction Criterion]
The drift \( w_{a,1}(n) \) measures the bit-length change of a symbolic block.
Averaging over all odd inputs yields:
\[
\mathbb{E}[w_{a,1}(n)] < 0 \quad \Longrightarrow \quad \text{global symbolic contraction}.
\]
This occurs if the trailing-zero decay outweighs the logarithmic growth from repeated \( T_3 \)-steps.
\end{remark}

\begin{figure}[H]
    \centering
    \includegraphics[width=0.9\textwidth]{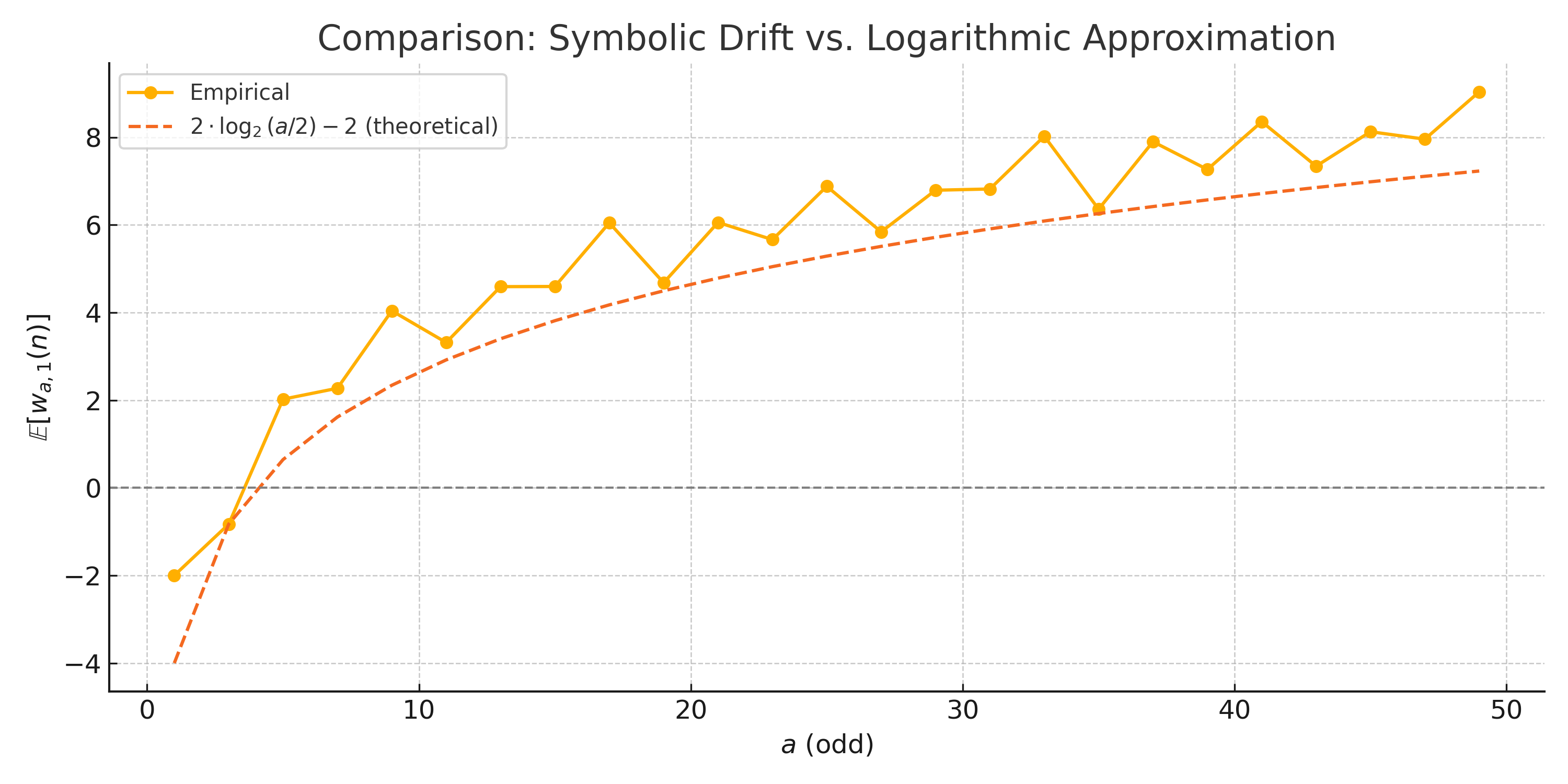}
    \caption{Empirical symbolic drift \( \mathbb{E}[w_{a,1}(n)] \) for odd \( a \leq 50 \), compared to the theoretical curve \( 2 \log_2(a/2) - 2 \). Minimum occurs at \( a = 3 \).}
    \label{fig:drift-a}
\end{figure}

\begin{remark}[Logarithmic Structure]
The shape of \( \mathbb{E}[w_{a,1}(n)] \) closely follows a logarithmic law:
\[
\mathbb{E}[w_{a,1}(n)] \approx 2 \cdot \log_2\left( \frac{a}{2} \right) - 2
\]
This arises because both \( \mathbb{E}[\mathrm{to}(n)] \approx 2 \) and \( \mathbb{E}[\mathrm{tz}(n')] \approx 2 \) are nearly constant.
\end{remark}

\begin{remark}[Symbolic Stability Principle]
Among all affine systems \( T_{a,1}(n) = \frac{a n + 1}{2^k} \), only the classical Collatz case \( a = 3 \) yields:
\[
\mathbb{E}[w_{3,1}(n)] \approx -0.83008 < 0
\]
This makes it the unique globally contracting system in the symbolic framework.
\end{remark}

\subsection{Symbolic Evaluation for \boldmath$n = 32$ (Even Case)}
\label{appendix:even-example}

We illustrate the symbolic digitwise emulation for the even input \( n = 32 \), where \( T(32) = 16 \).

The decimal expansion is:
\[
n = 2 \cdot 10^0 + 3 \cdot 10^1 \quad \Rightarrow \quad r_0 = 2,\quad r_1 = 3
\]

We now compute the updated digits from right to left:

\begin{align*}
p_0 &= r_1 \bmod 2 = 3 \bmod 2 = 1 \\
c_0 &= 0 \quad \text{(by definition)} \\
r_0' &= \delta_{\text{even}}(2, p_0, c_0) = \left\lfloor \frac{2 - 0}{2} \right\rfloor + 5 \cdot 1 = 1 + 5 = 6 \\[1ex]
p_1 &= 0 \quad \text{(no digit after \( r_1 \))} \\
c_1 &= r_1 \bmod 2 = 3 \bmod 2 = 1 \\
r_1' &= \delta_{\text{even}}(3, p_1, c_1) = \left\lfloor \frac{3 - 1}{2} \right\rfloor + 5 \cdot 0 = 1 + 0 = 1
\end{align*}

We obtain the updated digit sequence \( (r_0', r_1') = (6, 1) \), which corresponds to:
\[
T(32) = 6 \cdot 10^0 + 1 \cdot 10^1 = 16
\]

\begin{remark}
This confirms that the symbolic transition faithfully emulates the even-case arithmetic of the Collatz map at the digit level, including parity-aware carry handling.
\end{remark}

\section{Symbolic State Evolution Example}

The following table illustrates the symbolic digit evolution of a Collatz trajectory starting from \( n = 13 \).  
Each row corresponds to a timestep, showing the decimal value, and the symbolic states of the least and next-most significant digits.

\begin{table}[H]
\centering
\caption{Symbolic FSM states during the Collatz trajectory of \( n = 13 \)}
\begin{tabular}{|c|c|c|c|c|}
\hline
\textbf{Step} & \textbf{Value} & \textbf{State \( s_1 \)} & \textbf{State \( s_0 \)} \\
\hline
0    & 13  & (1,0,1) & (3,1,0) \\
1    & 40  & (4,0,0) & (0,0,0) \\
2    & 20  & (2,0,0) & (0,0,0) \\
3    & 10  & (1,0,0) & (0,1,0) \\
4    & 5   & (0,0,0) & (5,0,0) \\
5    & 16  & (1,0,0) & (6,1,0) \\
6    & 8   & (0,0,0) & (8,0,0) \\
7    & 4   & (0,0,0) & (4,0,0) \\
8    & 2   & (0,0,0) & (2,0,0) \\
9    & 1   & (0,0,0) & (1,0,0) \\
\hline
\end{tabular}
\label{tab:symbolic-fsm-example}
\end{table}

\begin{remark}
The symbolic configuration collapses progressively: from step 2 onward, all states in \( s_1 \) evolve toward the null state \((0,0,0)\), while \( s_0 \) ultimately enters the terminal cycle described in Theorem~\ref{thm:cycle}.
\end{remark}

\newpage

\section{Full Transition Table of the Symbolic FSM}

Table~\ref{tab:fsm} lists all possible state transitions for the symbolic digit automaton.
Each state is of the form $(r, p, c)$ and maps deterministically to its next state under either division or multiplication by $3n + 1$. The FSM is complete over the symbolic space $\mathcal{S} = \{0,\dots,9\} \times \{0,1\} \times \{0,1,2\}$.

\begin{table}[H]
\centering
\scriptsize
\renewcommand{\arraystretch}{0.9}
\begin{adjustbox}{max width=\textwidth, max totalheight=0.92\textheight}
\begin{tabularx}{\textwidth}{@{\hskip 2pt}l@{\hskip 6pt}X@{\hskip 6pt}X@{\hskip 6pt}X@{\hskip 6pt}X@{\hskip 6pt}X}
\toprule
State & Div target 1 & Div target 2 & Mul3 target 1 & Mul3 target 2 & Mul3 target 3 \\
\midrule
  0,0,0 & 0,0,0 & 0,1,0 &  0,0,0 &  1,0,0 &  2,0,0 \\
  0,0,1 &     - &     - &  1,0,0 &      - &      - \\
  0,0,2 &     - &     - &  2,0,0 &      - &      - \\
  0,1,0 & 5,0,0 & 5,1,0 &  0,1,0 &  1,1,0 &  2,1,0 \\
  0,1,1 &     - &     - &  1,1,0 &      - &      - \\
  0,1,2 &     - &     - &  2,1,0 &      - &      - \\
  1,0,0 & 0,0,0 & 0,1,0 &  3,0,0 &  4,0,0 &  5,0,0 \\
  1,0,1 & 0,0,0 & 0,1,0 &  4,0,0 &      - &      - \\
  1,0,2 &     - &     - &  5,0,0 &      - &      - \\
  1,1,0 & 5,0,0 & 5,1,0 &  3,1,0 &  4,1,0 &  5,1,0 \\
  1,1,1 & 5,0,0 & 5,1,0 &  4,1,0 &      - &      - \\
  1,1,2 &     - &     - &  5,1,0 &      - &      - \\
  2,0,0 & 1,0,0 & 1,1,0 &  6,0,0 &  7,0,0 &  8,0,0 \\
  2,0,1 &     - &     - &  7,0,0 &      - &      - \\
  2,0,2 &     - &     - &  8,0,0 &      - &      - \\
  2,1,0 & 6,0,0 & 6,1,0 &  6,1,0 &  7,1,0 &  8,1,0 \\
  2,1,1 &     - &     - &  7,1,0 &      - &      - \\
  2,1,2 &     - &     - &  8,1,0 &      - &      - \\
  3,0,0 & 1,0,0 & 1,1,0 &  0,1,0 &  1,1,0 &  9,0,0 \\
  3,0,1 & 1,0,0 & 1,1,0 &  0,1,0 &      - &      - \\
  3,0,2 &     - &     - &  1,1,0 &      - &      - \\
  3,1,0 & 6,0,0 & 6,1,0 &  0,0,0 &  1,0,0 &  9,1,0 \\
  3,1,1 & 6,0,0 & 6,1,0 &  0,0,0 &      - &      - \\
  3,1,2 &     - &     - &  1,0,0 &      - &      - \\
  4,0,0 & 2,0,0 & 2,1,0 &  2,1,0 &  3,1,0 &  4,1,0 \\
  4,0,1 &     - &     - &  3,1,0 &      - &      - \\
  4,0,2 &     - &     - &  4,1,0 &      - &      - \\
  4,1,0 & 7,0,0 & 7,1,0 &  2,0,0 &  3,0,0 &  4,0,0 \\
  4,1,1 &     - &     - &  3,0,0 &      - &      - \\
  4,1,2 &     - &     - &  4,0,0 &      - &      - \\
  5,0,0 & 2,0,0 & 2,1,0 &  5,1,0 &  6,1,0 &  7,1,0 \\
  5,0,1 & 2,0,0 & 2,1,0 &  6,1,0 &      - &      - \\
  5,0,2 &     - &     - &  7,1,0 &      - &      - \\
  5,1,0 & 7,0,0 & 7,1,0 &  5,0,0 &  6,0,0 &  7,0,0 \\
  5,1,1 & 7,0,0 & 7,1,0 &  6,0,0 &      - &      - \\
  5,1,2 &     - &     - &  7,0,0 &      - &      - \\
  6,0,0 & 3,0,0 & 3,1,0 &  0,0,0 &  8,1,0 &  9,1,0 \\
  6,0,1 &     - &     - &  9,1,0 &      - &      - \\
  6,0,2 &     - &     - &  0,0,0 &      - &      - \\
  6,1,0 & 8,0,0 & 8,1,0 &  0,1,0 &  8,0,0 &  9,0,0 \\
  6,1,1 &     - &     - &  9,0,0 &      - &      - \\
  6,1,2 &     - &     - &  0,1,0 &      - &      - \\
  7,0,0 & 3,0,0 & 3,1,0 &  1,0,0 &  2,0,0 &  3,0,0 \\
  7,0,1 & 3,0,0 & 3,1,0 &  2,0,0 &      - &      - \\
  7,0,2 &     - &     - &  3,0,0 &      - &      - \\
  7,1,0 & 8,0,0 & 8,1,0 &  1,1,0 &  2,1,0 &  3,1,0 \\
  7,1,1 & 8,0,0 & 8,1,0 &  2,1,0 &      - &      - \\
  7,1,2 &     - &     - &  3,1,0 &      - &      - \\
  8,0,0 & 4,0,0 & 4,1,0 &  4,0,0 &  5,0,0 &  6,0,0 \\
  8,0,1 &     - &     - &  5,0,0 &      - &      - \\
  8,0,2 &     - &     - &  6,0,0 &      - &      - \\
  8,1,0 & 9,0,0 & 9,1,0 &  4,1,0 &  5,1,0 &  6,1,0 \\
  8,1,1 &     - &     - &  5,1,0 &      - &      - \\
  8,1,2 &     - &     - &  6,1,0 &      - &      - \\
  9,0,0 & 4,0,0 & 4,1,0 &  7,0,0 &  8,0,0 &  9,0,0 \\
  9,0,1 & 4,0,0 & 4,1,0 &  8,0,0 &      - &      - \\
  9,0,2 &     - &     - &  9,0,0 &      - &      - \\
  9,1,0 & 9,0,0 & 9,1,0 &  7,1,0 &  8,1,0 &  9,1,0 \\
  9,1,1 & 9,0,0 & 9,1,0 &  8,1,0 &      - &      - \\
  9,1,2 &     - &     - &  9,1,0 &      - &      - \\
\bottomrule
\end{tabularx}
\end{adjustbox}
\caption{All valid state transitions in the symbolic FSM for all $60$ states $(r, p, c) \in \mathcal{S}$.}
\label{tab:fsm}
\end{table}

\section{Symbolic FSM Variants}
\label{app:binaryfsm}
\subsection{Quotient-Based FSM (Mod 2)}

The Mod-2 FSM encodes how the transformation $3n+1$ modifies the binary decay of $n$. Instead of computing values, it interprets Collatz dynamics structurally via symbolic rules over quotient sequences $q_i$ and $b_i$.

Let $q = (q_i)$ be the successive divisions of $n$ by 2, and $b = (b_i)$ the quotient chain of $3n+1$. For each index $i$, the FSM assigns a symbolic label from:
\[
\boxed{3n},\quad \boxed{3n + 1},\quad \boxed{3n + 2}
\]
using only local parity and division behavior. The rules are:

\begin{itemize}
  \item \textbf{Rule 1 (Initial step):} If $i = 0$, assign $\boxed{3n + 1}$.
  \item \textbf{Rule 2 (Parity mismatch):} If $\mathrm{par}(q_i) \ne \mathrm{par}(b_i)$, then $\boxed{3n + 1}$.
  \item \textbf{Rule 3 (Clean division, parity match):} If $q_i // 2 = q_{i+1}$:
  \[
  \begin{cases}
  b_i = 3q_i     &\Rightarrow \boxed{3n} \\
  b_i = 3q_i + 2 &\Rightarrow \boxed{3n + 2}
  \end{cases}
  \]
  \item \textbf{Rule 4 (Dirty division):} If $q_i // 2 \ne q_{i+1}$, assign $\boxed{3n + 2}$.
  \item \textbf{Rule 5 (Tail values):} For $b_i \in \{2,1\}$, use explicit matching.
\end{itemize}

This decoder assigns symbolic meanings to how $3n+1$ transforms the structure of $n$, without arithmetic simulation.

\subsubsection{Example: $n = 31$}

\begin{align*}
q &= [31,\ 15,\ 7,\ 3,\ 1] \\
b &= [94,\ 47,\ 23,\ 11,\ 5,\ 2,\ 1]
\end{align*}

\begin{center}
\renewcommand{\arraystretch}{1.2}
\begin{tabular}{c|c|c|c|c|l}
Step $i$ & $q_i$ & $b_i$ & $\mathrm{par}(q_i), \mathrm{par}(b_i)$ & $q_i // 2 = q_{i+1}$ & Operation \\\hline
0 & 31 & 94 & 1, 0 & --        & $3n + 1$ \\
1 & 15 & 47 & 1, 1 & \xmark    & $3n + 2$ \\
2 & 7  & 23 & 1, 1 & \xmark    & $3n + 2$ \\
3 & 3  & 11 & 1, 1 & \xmark    & $3n + 2$ \\
4 & 1  & 5  & 1, 1 & \xmark    & $3n + 2$ \\
5 & 0  & 2  & 0, 0 & --        & induced \\
6 & 0  & 1  & 0, 1 & --        & induced \\
\end{tabular}
\end{center}

\noindent\textbf{Symbolic sequence:} $\boxed{+1,\ +2,\ +2,\ +2,\ +2}$

The trailing entries $b_5 = 2$ and $b_6 = 1$ are not mapped from $q_i$, but are induced by the binary expansion of $3n+1$.

\subsubsection{Structural Use}

This symbolic decoder enables:

\begin{itemize}
  \item Structural classification of integers via FSM-derived symbol sequences
  \item Prediction of expansion/contraction via drift:
    \[
    w(n) = \mathrm{to}(n) \cdot \log_2\left(\frac{3}{2}\right) - \mathrm{tz}(n')
    \]
  \item Recognition of special forms (e.g., $2^k - 1$) with high symbolic lifespan
  \item Logical, rather than arithmetic, interpretation of the $3n+1$ process
\end{itemize}

Together with the digit-level FSM, the Mod-2 FSM provides a layered structural view of Collatz dynamics.

\subsection{Local Bit-Level FSM}

The bit-level FSM eliminates the need to compute $3n+1$. It operates directly on the binary form of $n$, using local rules and a sliding window.

\begin{itemize}
  \item \textbf{Input:} binary string $b_{k-1} \dots b_0$ of odd $n$
  \item \textbf{Step:} For each window $(b_{i+1}, b_i)$:
  \[
  \begin{aligned}
  00 &\mapsto 0 \\
  01 &\mapsto 0 \\
  10 &\mapsto 1 \\
  11 &\mapsto 1
  \end{aligned}
  \]
  \item \textbf{Growth prediction:}
    \[
    \Delta\mathrm{bits}(3n+1) =
    \begin{cases}
    +2 & \text{if } 3n+1 \geq 2^{k+1} \\
    +1 & \text{otherwise}
    \end{cases}
    \]
  \item \textbf{Symbol label:} 
    \[
    \texttt{10} \Rightarrow \boxed{+2},\quad \texttt{1} \Rightarrow \boxed{+1}
    \]
  \item \textbf{Result:} FSM yields symbolic string over $\{+2, +1, /2, \dots\}$, purely from $n$'s bits
\end{itemize}

\subsubsection{Example: $n = 161 = \texttt{10100001}_2$}

\begin{itemize}
  \item $3n+1 = 484 = \texttt{111100100}_2$, bit length grows from $8$ to $9$
  \item Symbol assigned: $\boxed{+1}$
  \item Followed by 3 halving steps: $\boxed{/2,\ /2,\ /2}$
\end{itemize}

\subsubsection*{Implication}

Each symbolic growth (e.g. $+2$) structurally creates trailing zeros → these must be consumed via division. Thus:

\[
\exists\ k\ \text{with } w(T^k(n)) < 0
\Rightarrow \text{Convergence}
\]

\subsubsection{Symbolic Length as Structural Predictor}

The bit-level FSM not only encodes symbolic structure, but also permits a \textbf{deterministic symbolic determination} of the exact number of steps until convergence. While not a closed-form formula, the symbolic step count
\[
L(n) := |\sigma(n)|
\]
can be computed by applying only local bit-level rules—without evaluating $3n + 1$ or simulating classical values.

At each stage, the number of trailing ones $\mathrm{to}(n)$ determines a symbolic growth phase, followed by trailing zeros $\mathrm{tz}(n')$ that encode a decay phase. Repeating this symbolic process yields:
\[
L(n) = \sum_{k=0}^{m} \left( \mathrm{to}(n_k) + \mathrm{tz}(n_k') \right),
\]
where $n_{k+1}$ is the result of structurally applying $\mathrm{to}(n_k)$ T3-steps and $\mathrm{tz}(n_k')$ halving steps.

\paragraph{Theorem.} \textit{For any odd $n$, the total symbolic step count $L(n)$ equals the number of symbolic FSM transitions from $n$ to $1$, and can be computed without evaluating $3n + 1$.}

This construction replaces arithmetic simulation with a fully symbolic recursion driven by bit patterns. While not an algebraic shortcut, it is a structurally grounded and mechanizable prediction mechanism.

\subsection{Comparison}

Both FSMs symbolically encode the $3n+1$ dynamics.

\begin{itemize}
  \item The \textbf{Quotient FSM} reconstructs the $3n+1$ action by comparing division chains $q_i$ and $b_i$
  \item The \textbf{Bit FSM} predicts symbolic drift directly from the binary representation of $n$
\end{itemize}

Both yield the same drift function $w(n)$ and confirm symbolic convergence. The bit-level FSM, however, requires no access to $3n+1$, making it truly predictive, symbolic, and deterministic.

\end{document}